\pdfoutput=1
\documentclass[3p,12pt]{elspreprint} 

\usepackage[utf8]{inputenc} 
\usepackage[english]{babel}  
\usepackage[T1]{fontenc}
\usepackage{lmodern}

\usepackage[intlimits]{amsmath}
\usepackage{amssymb} 
\usepackage{amsthm}
\usepackage{amsfonts}
\usepackage{bm}
\usepackage{mathtools} 

\usepackage{siunitx} 

\usepackage{setspace}
\usepackage{xspace}

\usepackage[noend]{algorithmic}
\usepackage{algorithm}

\usepackage{comment}
\usepackage{booktabs} 

\usepackage[labelformat=simple]{subcaption}   	

\usepackage{hyperref} 

\newcommand{\Lame}{Lam{\'e}\xspace} 

\newcommand{\myMat}[1]{\mathbf{#1}}

\newcommand{\myVec}[1]{\mathbf{#1}}

\newcommand{\myNorm}[1]{\lVert #1 \rVert}
\newcommand{\myCardinality}[1]{\# #1}

\newcommand\indexA{i} 		
\newcommand\totalA{I} 		
\newcommand\indexB{j} 		
\newcommand\totalB{J} 		
\newcommand\indexC{k} 		
\newcommand\indexD{m} 		
\newcommand\indexE{n} 		
\newcommand\indexG{\ell} 		
\newcommand\totalG{L} 		

\newcommand\indexSpan{s}	
\newcommand\indexSpanTrimmed{t}	

\newcommand{\Bspline}{B} 		
\newcommand{\BsplineSeg}{\mathcal{B}} 	
\newcommand{\KV}{\varXi} 		
\newcommand{\KVRefined}{\hat{\varXi}} 	
\newcommand\uu{r} 			
\newcommand\uuRefined{\hat{\uu}} 	
\newcommand\vv{\uu} 			
\newcommand\UVsurf{\pt{\uu}} 	
\newcommand\uusurf{\uu_1} 		
\newcommand\vvsurf{\vv_2} 		
\newcommand\multi{m} 			
\newcommand{\knotinsertion}{d_k}

\newcommand{\CP}{\pt{c}}		
\newcommand{\CPRefined}{\hat{\CP}} 	
\newcommand{\CPara}{c}			
\newcommand{\CPCoons}{\CP^{\uu}}	
\newcommand\collocationPoint{\pt{x}^c}	

\newcommand\abscissa{\bar{\uu}}			
\newcommand\greville{\abscissa}		

\newcommand{\BsplineSpace}{\mathbb{S}} 

\newcommand\pu{p}			
\newcommand\pusurf{\pu_1} 		
\newcommand\pvsurf{\pu_2} 		

\newcommand{\curve}{{\tau}}
\newcommand{\curveRefined}{\hat{{\tau}}}
\newcommand{\curveBez}{{\tau}_{b}}
\newcommand{\curveBezOpposite}{{\tau}_{e}}

\newcommand\ieKs{\tau}                  
\newcommand\ieLocal{\tilde{\ieKs}} 	
\newcommand\ieReference{\grave{\ieKs}}  

\newcommand{\exBspline}{B^e} 		
\newcommand{\degSet}{\mathbb{J}} 	
\newcommand{\exteriorSet}{\mathbb{K}} 	
\newcommand{\NSet}{\mathbb{T}} 	  	
\newcommand{\eweight}{e} 	  	
\newcommand{\Nbasis}{\psi} 	  	
\newcommand{\taylorpoint}{\tilde{\uu}} 
\newcommand{\poly}{\tilde{\alpha}} 		
\newcommand{\lambdapoint}{\mu} 	
\newcommand{\trim}{t}
\newcommand{\variableKnot}{\hat{\uu}}

\newcommand{\PrimaryBasis}{\varphi} 	
\newcommand{\DualBasis}{\psi} 		

\newcommand{\DOF}{n} 			
\newcommand{\domain}{\Omega}
\newcommand{\patchdomain}{\mathcal{A}^{\textnormal{v}}}
\newcommand{\boundary}{\Gamma}
\newcommand{\Ltwo}{L_{2}}

\newcommand\primary{u} 			
\newcommand\dual{t} 			
\newcommand\utens{\primary} 		
\newcommand\ttens{\dual} 		

\newcommand{\R}{\mathbb{R}}

\DeclareMathOperator{\supp}{supp}

\newcommand\operate[1]{(#1)}
\newcommand\dx[1]{\;\text{d}#1}

\providecommand\url[1]{\emph{#1}}

\newcommand\mymax[1]{\max\{#1\}}
\newcommand\fromto[2]{\{ #1, \dots, #2 \}}

\newcommand\dgamma[1]{\dx\mathrm{s}_{\pt{#1}}}

\newcommand\tens[2]{\mathsf{#1}}
\newcommand{\trans}{\intercal} 

\newcommand\fund[1]{\tens{#1}{2}}

\newcommand\beistrich{,\,}
\newcommand\und{\,\text{and}\,}

\newcommand\with{\,\text{with}\,}
\newcommand\for{\,\text{for}\,}

\newcommand\pt[1]{\boldsymbol{#1}}
\newcommand\ofpt[1]{(\pt{#1})}
\newcommand\sourcept{\tilde{\pt{x}}}

\newcommand\boundarypt{\pt{y}}

\newcommand\err{\epsilon}

\newcommand\mychi{\mathcal{X}}

\newcommand{\tikzfig}[5]{
  \begin{figure}[thb!]
    \centering
    \includegraphics{#5}
    \caption{#3}
    \label{#4}
  \end{figure}
}

\newcommand{\tikzfigposition}[6]{
  \begin{figure}[#6]
    \centering
    \includegraphics{#5}
    \caption{#3}
    \label{#4}
  \end{figure}
}

\newcommand{\tikzsubfig}[5]{

    \centering
    \includegraphics{#5}%
    \subcaption{#3}
    \label{#4}
}

\newcommand{\tikzsubfigscale}[5]{

    \centering
    \includegraphics[scale=#2]{#5}%
    \subcaption{#3}
    \label{#4}
}

{
  \begin{table}[htb]
    \def\mycap{#1}
    \def\mylabel{tab:#2}
    \centering
    \begin{tabular}{#3}
      \toprule
}
{
  \bottomrule
  \end{tabular}
  \caption{\mycap}
  \label{\mylabel}
  \end{table}
}
\newcommand{\mytableheader}[1]{
       #1 \\ \midrule
}

\newenvironment{mytableposition}[4]%
{
  \begin{table}[#4]
    \def\mycap{#1}
    \def\mylabel{tab:#2}
    \centering
    \begin{tabular}{#3}
      \toprule
}
{
  \bottomrule
  \end{tabular}
  \caption{\mycap}
  \label{\mylabel}
  \end{table}
}

\newcommand{\store}{=} 	
 
\newenvironment{myalgorithm}[2]%
{
  \begin{algorithm}
    \caption{#1}
    \label{#2}
    \begin{algorithmic}[1]
}   
{
\end{algorithmic}
\end{algorithm}
}

\ifdefined\TodoNotes
    \newcommand{\mycomment}[2][-] 
    {
        {%
            \ifdefined\tikzExtOn
            \tikzset{external/export=false}
            \fi
            \setstretch{0.1}
            \todo[color={green!33}]{ \textbf{[\uppercase{#1}]} #2}
            \pdfcomment[color=yellow,opacity=0.0,author=#1]{}
        }%
    }	
    \newcommand{\myhighlight}[2][-] 
    {
    \pdfcomment[color=yellow,voffset=-5pt,hoffset=3pt,opacity=0.0,author=#1]{}
    \hl{#2}%
    }	
    \setstcolor{red} 
    \newcommand{\mydelete}[2][-] 
    {
    \pdfcomment[color=yellow,voffset=-5pt,hoffset=3pt,opacity=0.0,author=#1]{}
    \st{#2}%
    }
\else
    \newcommand{\mycomment}[2][-]{}
    \newcommand{\myhighlight}[2][-]{}
    \newcommand{\mydelete}[2][-]{}
\fi

\ifdefined\tikzSubmission
\newcommand\revdel[2]{}
\newcommand\revadd[2]{}
\newcommand\revmod[2]{}
\newcommand\revcomment[1]{} 
\else
\newcommand\revdel[2]{ %
{
  \pdfmarkupcomment[markup=StrikeOut,color=red,author=Marussig and Zechner]{#1}{#2}
}
}
\newcommand\revadd[2]{ %
{
  \pdfmarkupcomment[markup=Underline,color=green,author=Marussig and Zechner]{#1}{#2}
}
}
\newcommand\revmod[2]{ %
{
  \pdfmarkupcomment[markup=Highlight,color=yellow,author=Marussig and Zechner]{#1}{#2}
}
}
\newcommand\revcomment[1]{ %
{
  \pdfcomment[color=yellow,author=Marussig and Zechner,voffset=8pt,opacity=0.8]{#1}
}
}
\fi

\newtheoremstyle{myremark}
{3pt}
{3pt}
{}
{}
{\itshape}
{:}
{.5em}
{}
\theoremstyle{myremark}

\newcommand{\myEqref}[1]{Equation~(\ref{#1})}
\newcommand{\myeqref}[1]{equation~(\ref{#1})}	
\newcommand{\myfigref}[1]{Figure~\ref{#1}}
\newcommand{\myFigref}[1]{Figure~\ref{#1}}
\newcommand{\myRef}[1]{(\ref{#1})}
\newcommand{\mytabref}[1]{Table~\ref{#1}}

\newcommand{\mySecref}[1]{Section~\ref{#1}}
\newcommand{\mysecref}[1]{Section~\ref{#1}}
\newcommand{\myalgref}[1]{Algorithm~\ref{#1}}
\newcommand{\myappref}[1]{\ref{#1}}

\begin{document}

%
  
\title{Stable Isogeometric Analysis of Trimmed Geometries}

\begin{frontmatter}

\author[ifbaddr]{Benjamin Marussig\corref{cor1}}
\author[ifbaddr]{Jürgen Zechner}
\author[ifbaddr,newcastleaddr]{Gernot Beer}
\author[ifbaddr]{Thomas-Peter Fries}

\address[ifbaddr]{Institute of Structural Analysis, Graz University
  of Technology, Lessingstraße 25/II, 8010 Graz, Austria}

\address[newcastleaddr]{Centre for Geotechnical and Materials Modelling, University of Newcastle,
  Callaghan, NSW 2308, Australia}


\cortext[cor1]{Corresponding author. 
  Tel.: +43 316 873 6181, fax: +43 316 873 6185, mail: \url{ifb@tugraz.at}, web: \url{www.ifb.tugraz.at}}

\begin{abstract}
    We explore extended B-splines as a stable basis for isogeometric analysis with trimmed parameter spaces.
    The stabilization is accomplished by an appropriate substitution of B-splines that may lead to ill-conditioned system matrices.
    The construction for non-uniform knot vectors is presented.
    The properties of extended \mbox{B-splines} are examined in the context of interpolation, potential, and linear elasticity problems and excellent results are attained. 
    The analysis is performed by an isogeometric boundary element formulation using collocation.
    It is argued that extended B-splines provide a flexible and simple stabilization scheme which ideally suits the isogeometric paradigm.
\end{abstract}


\end{frontmatter}


%
  \section{Introduction}
\label{sec:introduction}

Trimmed tensor product surfaces are often used in Computer Aided Geometric Design (CAGD) models because they offer a convenient way to efficiently represent
non-rectangular surface topologies.
However, their application to an analysis requires further consideration.
First of all, trimming procedures are used to define visible areas~$\patchdomain$ over surfaces independent of the underlying parameter space.
While the parameter space is determined by a regular grid of knot spans, the actual shape of~$\patchdomain$ is arbitrarily defined by trimming curves.
An example of a trimmed surface is shown in \myfigref{fig:trimmedPatch}.
It should be noted that the mathematical description, i.e.~the tensor product basis \ref{fig:trimmedBasis} and the related control grid of the original surface \ref{fig:tensorProductPatch}, does not change.
Hence, it is not appropriate to define integration elements based on the knot spans of the parameter space, as it is usually done in isogeometric analysis.
\begin{figure}[htb]
  \centering
  \begin{minipage}[b]{0.31\textwidth}
      \centering
      \tikzsubfig{tikz/NURBSpatch.tex}{0.20}{Original Surface}{fig:tensorProductPatch}{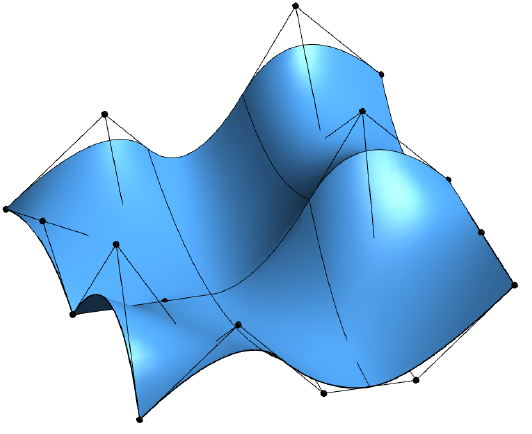}
  \end{minipage}
  \begin{minipage}[b]{0.31\textwidth}
      \centering
      \tikzsubfig{tikz/NURBSpatchtrimmedbasis.tex}{1.2}{Trimmed Parameter Space}{fig:trimmedBasis}{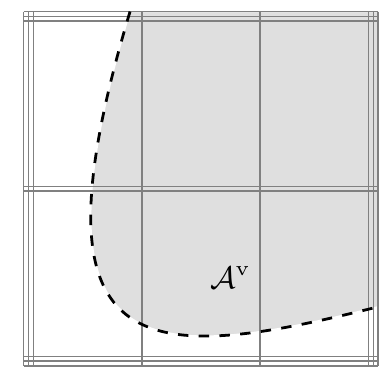}
  \end{minipage}
  \begin{minipage}[b]{0.31\textwidth}
      \centering
      \tikzsubfig{tikz/NURBSpatchtrimmed.tex}{0.20}{Trimmed Surface}{fig:trimmedPatch}{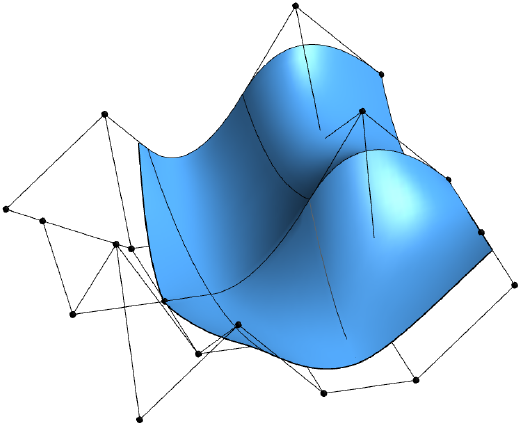}
  \end{minipage}
  \caption{Trimmed tensor product surface: (a) regular surface defined by a tensor product basis, (b) trimmed parameter space where a trimming curve (dashed line) specifies the visible part $\patchdomain$ of (c) the resulting trimmed surfaces.}
  \label{fig:tensorProductPatchTrimmed}
\end{figure}
Moreover, a trimmed basis contains \emph{degenerated} basis functions which exist only partially within~$\patchdomain$.
Their supports may be very small and thus the condition number of the resulting system matrices can become arbitrarily large.
Consequently, a trimmed basis is not guaranteed to be stable~\cite{Hoellig2012a}.

This stability aspect has scarcely been considered in previous works related to the analysis of trimmed geometries.
Current approaches focus mainly on proper integration of trimmed CAGD models and may be summarized as global and local techniques.
The former \cite{Beer2015a,Harbrecht2010a} substitutes a trimmed surface by one or several regular ones.
As a consequence, the continuity is either locally reduced by the non-smooth edges introduced or the regular surfaces are distorted in many topological cases.
Moreover, this approach cannot be applied to real-world geometries without user intervention.
The other techniques use the parameter space of the trimmed object as background parametrization and the trimming curve determines which part is considered for the analysis, i.e.~$\patchdomain$.
All this information is provided by the CAGD model, but the analysis of knot spans cut by the trimming curve requires special attention.
There are different approaches \cite{Kim2009a,Kim2010a,Schmidt2012a,Wang2013a} which locally substitute the trimmed area by regular elements providing a mapping from the reference element where quadrature points are specified to the trimmed parameter space. 
In addition, tailored integration formulae may be used \cite{Nagy2014a,wang2015a}.
Due to their local nature, these concepts can be applied to complex CAGD models.
However, the resulting basis contains degenerated basis functions that are truncated at the trimming curve.
Besides the resulting stability aspect, this is a disadvantage in the context of collocation methods, because collocation points of degenerated basis functions may be located outside of the actual domain.

In this paper, it is proposed to address the issues of local schemes by the application of \emph{extended B-splines}. 
These splines have been introduced in the context of weighted extended B-splines (WEB-splines)~\cite{Hoellig2003b,Hoellig2003a,Hoellig2001a,Hoellig2002a}, a fictitious domain method based on B-splines.
The term \emph{weighted} indicates that a weighting function is used to define the domain of interest within a regular B-spline basis. 
In addition, it allows the enforcing of essential boundary conditions along the boundary of the domain. 
Similar to trimmed geometries, the resulting basis contains degenerated basis functions and requires a stabilization, which is accomplished by the concept of extended B-splines.
Although the application of WEB-splines to trimmed geometries has already been considered in the landmark paper of isogeometric analysis~\cite{Hughes2005a}, only a few attempts were made to unite these fields:
a combination of isogeometric analysis and WEB-splines is proposed in \cite{Hoellig2012a} focusing on the application of the weighting functions rather than the stabilization and recently, WEB-splines have been applied to isogeometric collocation with uniform parameter spaces \cite{Apprich2016a}.
Usually, extended B-spline applications are tailored to the uniform case, e.g.~\cite{Apaydin2008a,Rueberg2012a,SrinivasKumar2006a,Verschaeve2015aPre}, most likely because
this simplifies the construction. 
However, the concept is by no means restricted to uniform parameter spaces.
The mathematical foundation of the generalized non-uniform situation has been presented in \cite{Hoellig2003a} and in this work, we will address the computational aspect of setting up non-uniform extended B-splines.
In contrast to the previously mentioned attempts on combining extended B-splines and isogeometric analysis, this paper focuses on the stabilization ability of extended B-splines and their application to a collocation scheme within parameter spaces defined by non-uniform knot spacing and multiple knots.

The analysis of trimmed geometries is performed by means of a collocated isogeometric boundary element method (BEM). 
In the present context, BEM has the following advantages: (i) it does not need a volume discretization since it is based on a boundary representation just like CAGD models, (ii) it does not necessarily require interelement continuity \cite{Patterson1984a}, hence no coupling strategy for adjacent trimmed patches is needed, and (iii) all boundary conditions are generally employed in an integral sense.  
The latter argument, together with the fact that $\patchdomain$ is already defined by a trimming curve, permits to neglect the weighting function used in the WEB-spline approach. 
Hence, emphasis can be placed on the impact of extended B-splines.
To the best of the authors' knowledge, this is the first application of extended B-splines to a BEM formulation.
However, this paper is mainly dedicated to the stabilization of trimmed geometries and hence the BEM is only briefly reviewed.
For an in-depth discussion of the applied BEM implementation the interested reader is referred to \cite{marussig2016b,Marussig2014aCMAME}.

The paper begins by reviewing essential aspects of conventional B-splines in \mysecref{sec:NURBS}.
\mysecref{sec:ExtendedBsplines} goes on to discuss the general construction of extended B-splines from an engineering point of view.
\mySecref{sec:analysis} is dedicated to the actual analysis of trimmed geometries.
The performance of the proposed stabilization is critically assessed by various interpolation, potential, and elasticity problems in \mysecref{sec:results}. 
The paper closes with concluding remarks in \mysecref{sec:conclusion}. Supplementary information for setting up non-uniform extended B-splines is provided in \myappref{chap:extrapolationWeightsExample}.

  \section{B-splines}
\label{sec:NURBS}

This section gives a brief overview of B-splines focusing on aspects which are crucial for the subsequent derivation of extended B-splines and their application to analysis.
For further information related to spline theory, the interested reader is referred to \cite{boor2001b,Farin2002b}. 
Detailed descriptions of efficient algorithms can be found in~\cite{Piegl1997b}. 
In this paper, the term B-spline refers to basis functions, while the  objects obtained, i.e.~curves or surfaces, are generally denoted as patches.

\subsection{Basis Functions}
\label{sec:NURBSbasis}

\begin{figure}[tb!]
  \centering
  \begin{minipage}[b]{0.45\textwidth}
	\tikzsubfig{tikz/bsplineSegments.tex}{1.2}{$\KV=\left\{1,2,3,4\right\}$}
	{fig:BsplineSegment_ex1}{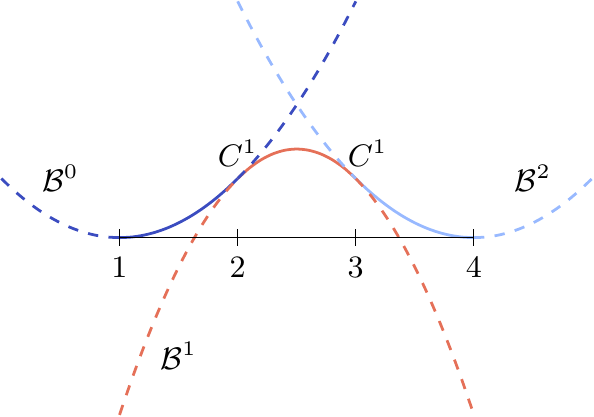}
  \end{minipage}
  \hfill
  \begin{minipage}[b]{0.45\textwidth}
	\tikzsubfig{tikz/bsplineSegments2.tex}{1.2}{$\KV=\left\{1,2.5,2.5,4\right\}$}
	{fig:BsplineSegment_ex2}{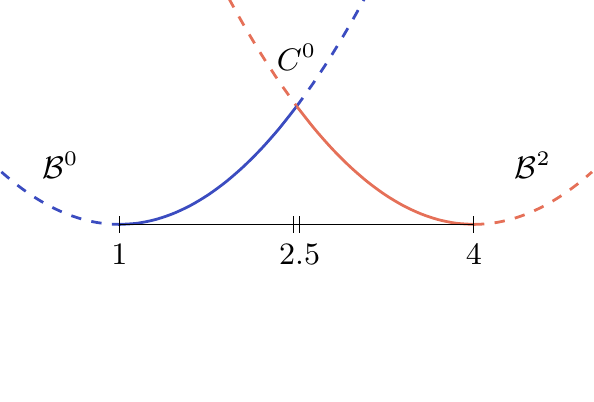}
  \end{minipage}
  \caption{Polynomial segments $\BsplineSeg^{\indexSpan}$ of a quadratic B-spline due to different knot vectors~$\KV$. The resulting polynomial segments are indicated by dashed lines, whereas solid lines represent the corresponding basis function.}
  \label{fig:BsplineSegment}
\end{figure}

B-splines $\Bspline_{\indexA,\pu}$ consist of polynomial segments which possess a certain smoothness at their connection.
For a fixed polynomial degree~$\pu$, they are defined recursively 
by a strictly convex combination of B-splines of the previous degree $\pu-1$ given by
\begin{align}
	\label{eq:Bspline_Np}
	\Bspline_{\indexA,\pu}(\uu) & = \frac{\uu-\uu_{\indexA}}{\uu_{\indexA+\pu}-\uu_{\indexA}} \: \Bspline_{\indexA,\pu-1}(\uu) 
	+ \frac{\uu_{\indexA+p+1}-\uu}{\uu_{\indexA+p+1}-\uu_{\indexA+1}} \: \Bspline_{\indexA+1,\pu-1}(\uu)
  \shortintertext{\with}
  \label{eq:Bspline_N0}
  \Bspline_{\indexA,0}(\uu) & = 
  \begin{cases} 
    1 & \textnormal{if } \uu_{\indexA}\leqslant \uu < \uu_{\indexA+1}\\
    0 & \textnormal{otherwise.}
  \end{cases}
\end{align}
The essential element for this construction is the \emph{knot vector}~$\KV$ characterized as a non-decreasing sequence of coordinates~{$\uu_\indexA \leqslant \uu_{\indexA+1}$}.
The parameters~$\uu_{\indexA}$ are termed \emph{knots} and the half-open interval $\left[\uu_{\indexA}, \uu_{\indexA+1}\right)$ is called $\indexA$th \emph{knot span}.
Each knot span has $\pu+1$ non-vanishing B-splines.
Each basis function is entirely defined by $\pu + 2$ knots and its support, $\supp{ \{\Bspline_{\indexA,\pu} \} }=\fromto{\uu_{\indexA}}{\uu_{\indexA+\pu+1}}$, is local.
Within each non-zero knot span~$\indexSpan$, $\uu_{\indexSpan} < \uu_{\indexSpan+1}$, of its support,  $\Bspline_{\indexA,\pu}$ is described by a polynomial segment $\BsplineSeg^{\indexSpan}_{\indexA}$.
Each knot value indicates a location within the parameter space which is not $C^{\infty}$-continuous, i.e.~where two adjacent $\BsplineSeg^{\indexSpan}_{\indexA}$ join.
Successive knots may share the same value, which is indicated by the knot multiplicity~$\multi$,  i.e.~$\uu_{\indexA} = \uu_{\indexA+1} = \dots = \uu_{\indexA+\multi-1}$.
In general, the continuity between adjacent segments is $C^{\pu-\multi}$.
The control of continuity is demonstrated for a quadratic B-spline in \myfigref{fig:BsplineSegment}.
If the multiplicity of the first and last knot is equal to $\pu+1$, the knot vector is denoted as \emph{open} knot vector.

A B-spline basis defined by a knot vector $\KV$ forms
a \emph{partition of unity}, i.e.
\begin{align}
	\label{eq:Bspline_partitionOfUnity}
	\sum_{\indexA=0}^{\totalA-1}  \Bspline_{\indexA,\pu}(\uu) = 1, && \uu \in \left[\uu_0,\uu_{\totalA+\pu}  \right]
\end{align}
and the basis functions are \emph{linear independent}, i.e. 
\begin{align}
	\label{eq:Bspline_linearIndependent}
	\sum_{\indexA=0}^{\totalA-1} \Bspline_{\indexA,\pu}(\uu) \: \CPara_\indexA  &= 0  
\end{align}
is satisfied if and only if $\CPara_\indexA = 0,~\indexA = 0,\dots,\totalA-1$.
Due to the latter property, every piecewise polynomial~$f_{\pu,\KV}$ of degree $\pu$ over a knot sequence~$\KV$ can be uniquely described by a linear combination of corresponding B-splines $\Bspline_{\indexA,\pu}$. 
Hence, they form a \emph{basis} of the space $\BsplineSpace_{\pu,\KV}$ collecting all such functions
\begin{align}
	\label{eq:splineSpace}
	\BsplineSpace_{\pu,\KV}  &= \sum_{\indexA=0}^{\totalA-1}  \Bspline_{\indexA,\pu} \: \CPara_\indexA, & \CPara_\indexA \in \R.
\end{align}

The first derivative of a B-spline is computed by a linear combination of B-splines of the previous degree 
\begin{align}
	\label{eq:Bspline_Np_Deriv}
	\Bspline^\prime_{\indexA,\pu}(\uu)  = \frac{\pu}{\uu_{\indexA+\pu}-\uu_{\indexA}} \: \Bspline_{\indexA,\pu-1}(\uu) 
				    - \frac{\pu}{\uu_{\indexA+\pu+1}-\uu_{\indexA+1}} \: \Bspline_{\indexA+1,\pu-1}(\uu).
\end{align}
For the computation of the $k$th derivative, this is generalized to
\begin{align}
	\label{eq:Bspline_Np_kDeriv}
	\Bspline^{(k)}_{\indexA,\pu}(\uu) &= 
	\frac{\pu!}{\left(\pu - k \right)! } 
	\sum_{\indexG=0}^{k}  
	a_{k,\indexG} \:
	\Bspline_{\indexA+\indexG,\pu-k}(\uu) 
	\shortintertext{\with} 
        \begin{split}
	a_{0,0} &= 1,  \\
	a_{k,0} &= \frac{a_{k-1,0}}{\uu_{\indexA+\pu-k+1}-\uu_{\indexA}},  \\
	a_{k,\indexG} &= \frac{a_{k-1,\indexG} - a_{k-1,\indexG-1}}{\uu_{\indexA+\pu+\indexG-k+1}-\uu_{\indexA+\indexG}}
	\qquad \indexG = 1,\dots,k-1 ,  \\
	a_{k,k} &= \frac{-a_{k-1,k-1}}{\uu_{\indexA+\pu+1}-\uu_{\indexA+k}}  .
        \end{split}
\end{align}

\subsection{Patches}
\label{sec:patches}

B-spline curves of degree $\pu$ are defined by basis functions $\Bspline_{\indexA,\pu}$ based on a knot vector~$\KV$ with corresponding coefficients in physical space $\CP_{\indexA}$ which are referred to as \emph{control points}.
The geometrical mapping~$\mychi\operate{\uu}$ from parameter space to physical space is given by
\begin{align}
  \label{eq:Bspline_mapping}
  \mychi\operate{\uu} &\coloneqq \pt{x}\operate{\uu} = 
  \sum_{\indexA=0}^{\totalA-1} \Bspline_{\indexA,\pu}\operate{\uu} \: \CP_{\indexA} 
\end{align}
with $\totalA$ representing the total number of basis functions.
In general, the control points $\CP_{\indexA}$ are not interpolatory, i.e.~they do not lie on the curve.
The related Jacobi-matrix is
\begin{align}
  \label{eq:Bspline_mapping_jac}
  \myMat{J}_\mychi \operate{\uu} &\coloneqq %
  \sum_{\indexA=0}^{\totalA-1} \Bspline^\prime_{\indexA,\pu}\operate{\uu} \: \CP_{\indexA} .
\end{align}

This concept can be easily extended to surfaces by means of tensor products.
Basis functions for B-spline surfaces are obtained by the tensor product of univariate B-splines which are defined by separate knot vectors $\KV_\totalA$ and $\KV_\totalB$.
These knot vectors determine the parametrization in the directions $\uusurf$ and $\vvsurf$, respectively. 
Moreover, they span the bivariate basis of a surface and specify its local coordinates $\pt{\uu}=(\uusurf,\vvsurf)^\trans$.
Combined with a bidirectional grid of control points $\CP_{\indexA,\indexB}$ the geometrical mapping is given by 
\begin{align}
	\label{eq:Bspline_mapping_patch}
	\mychi\operate{\UVsurf} &\coloneqq \pt{x}\operate{\uusurf,\vvsurf} = 
	\sum_{\indexA=0}^{\totalA-1} \sum_{\indexB=0}^{\totalB-1} \Bspline_{\indexA,\pusurf} \operate{\uusurf}  \Bspline_{\indexB,\pvsurf}\operate{\vvsurf}    \: \CP_{\indexA,\indexB} .
\end{align}
The polynomial degrees for each parametric direction are denoted by $\pusurf$ and $\pvsurf$.
The Jacobian of the mapping~\myRef{eq:Bspline_mapping_patch} is computed by substituting the occurring univariate \mbox{B-splines} by their first derivatives, alternately for each direction.
In general, derivatives of B-spline surfaces are specified by
\begin{align}
	\label{eq:Bspline_deriv_uu_patch}
	\frac{\partial^{k+l} }{ \partial^k \uusurf \partial^l \vvsurf}  \pt{x}\operate{\uusurf,\vvsurf} &=	
	\sum_{\indexA=0}^{\totalA-1} \sum_{\indexB=0}^{\totalB-1} \Bspline^{(k)}_{\indexA,\pusurf} \operate{\uusurf}  \Bspline^{(l)}_{\indexB,\pvsurf}\operate{\vvsurf}    \: \CP_{\indexA,\indexB}.
\end{align}

In general, a given patch $\curve$ can be refined such that the resulting object $\curveRefined$ is equivalent to the original one, i.e.~$\curve \equiv \curveRefined$.
The related procedures are called \emph{knot insertion} and \emph{degree elevation}.
In both cases, a given knot vector $\KV$ is extended to a knot vector $\KVRefined$ by adding new knots $\uuRefined$.
Hence, the number of basis functions is increased.
Moreover, the positions $\uuRefined$ determine the control points $\CPRefined_\indexA$ of the refined patch.
In other words, the refinement procedures define a new basis and set of control point coordinates without changing the geometry.
A comprehensive discussion of these refinement schemes can be found in \cite{Piegl1997b}.

\subsection{Spline Interpolation}

In case of a spline interpolation problem, a given function $f$ shall be approximated by a B-spline curve~$I_h f \coloneqq \sum_{\indexA=0}^{\totalA-1} { \Bspline_{\indexA,\pu} } \: \CPara_\indexA$.
They agree at $\totalA$ interpolation points~$\abscissa_\indexB$ if and only if
\begin{align}
	\label{eq:splineInterpolation}
	f(\abscissa_\indexB)  &= \sum_{\indexA=0}^{\totalA-1} { \Bspline_{\indexA,\pu}(\abscissa_\indexB) } \: \CPara_\indexA, & \indexB = 0,\dots,\totalA-1 &.
\end{align}
The corresponding system of equations consists of the unknown coefficients~$\CPara_\indexA$ and the \emph{spline collocation matrix}~$\myMat{A}_\uu$ which is defined by
\begin{align}
  \label{eq:splineCollocationMatrix}
  \myMat{A}_\uu  [\indexB,\indexA] &=  \Bspline_{\indexA,\pu}(\abscissa_\indexB), && \indexA,\indexB = 0,\dots,\totalA-1.
\end{align}
The \emph{Schoenberg-Whitney theorem} \cite{boor2001b,Farin2002b} states that the matrix $\myMat{A}_\uu$ is invertible if and only if 
\begin{align}
	\label{eq:schoenberg-whitney}
	\Bspline_{\indexA,\pu}(\abscissa_\indexA)& \neq 0, & \indexA = 0,\dots,\totalA-1 .
\end{align}
Thus, each interpolation point must be located within the support of its corresponding B-spline.
In general, the interpolation error for every continuous function~$f$ over a fixed interval~$[a,b]$ is bounded by
\begin{align}
	\label{eq:splineInterpolationError}
	\myNorm{ f - I_h f } &  
        \leq  C_{\pu} 
        \left( 1 + \myNorm{ I_h } \right)  \myNorm{ f^{(\pu + 1)} }_{\infty} \: \lvert \uu \rvert^{\pu + 1}  
	\shortintertext{\with}  
	\lvert \uu \rvert & \coloneqq \max_i \Delta\uu_i = \max_i \left( \uu_{i+1} - \uu_{i} \right) \\
	\myNorm{ I_h } & \coloneqq \mymax{ \myNorm{I_h f} / \myNorm{f} : f \in C[a,b] \backslash \{0 \} \: }  .
\end{align}
Proofs and more detailed information can be found in~\cite{boor2001b}. 
The factor $\lvert \uu \rvert$ of the bound~\myRef{eq:splineInterpolationError} indicates that the knot placement influences the approximation quality.
However, in order to find optimal knots, information about the target function~$f$ must be given, which is generally not the case in analysis. Hence, it is focused on the norm~$\myNorm{ I_h }$.

Since condition~\myRef{eq:schoenberg-whitney} guarantees that $\myMat{A}_\uu$ does not become singular, it is expected that $\myNorm{ I_h }$ gets large if a $\abscissa_\indexB$ approaches the limits of its allowed range.
Non-uniformity of the points $\abscissa_\indexB$ is another reason for an increasing $\myNorm{ I_h }$. 
In fact, $\myNorm{ I_h }$ gets arbitrary large if two interpolation points approach each other, while the others are fixed. 
Several authors \cite{Auricchia2010a,boor2001b,Li2011a} recommend interpolating at the \emph{Greville abscissae}, which are obtained by the knot average 
\begin{align}
    \label{eq:greville}
    \greville_{\indexA} & = \frac{\uu_{\indexA+1}+\uu_{\indexA+2} + \dots +\uu_{\indexA+\pu}}{\pu}.
\end{align}
This abscissae are well known in CAGD and used for different purposes, e.g.~to generate a linear geometrical mapping \cite{Farin2002b}.
The most important feature of this approach is that it induces a stable interpolation scheme for moderate degrees~$\pu$.
The only abscissae that provide a stable interpolation for any degree are proposed by Demko~\cite{Demko1985a}, but their computation is more involved than using the knot average~\myRef{eq:greville}.

Despite the preferred choice, the abscissae $\abscissa$ are generally denoted as \emph{anchors} for the remainder of the paper.
In particular, they are used as a means of linking basis functions $\Bspline_{\indexA,\pu}$ to a point at a specific location $\abscissa_{\indexA}$ in the parameter space.

\subsection{Quasi Interpolation}
\label{sec:splineInterpolation}

Quasi interpolation methods allow the computation of spline approximations without solving a linear system of equations.
Here the so-called \emph{de Boor--Fix} or dual functional $\lambda_{\indexB,\pu}$ \cite{boor2001b,deBoor1973aa} is introduced: 
for any piecewise polynomial $f \in \BsplineSpace_{\pu,\KV}$,
\begin{align}
	\label{eq:deBoorFix}
	f &= \sum^{\totalB-1}_{\indexB=0} \lambda_{\indexB,\pu}(f) \: \Bspline_{\indexB,\pu}
	\shortintertext{\with}  
	\label{eq:deBoorFixFunctional}
	\lambda_{\indexB,\pu}(f) &= \frac{1}{\pu!} \sum^\pu_{\indexC=0} (-1)^\indexC \Nbasis^{\left(\pu-\indexC\right)}_{\indexB,\pu}(\lambdapoint_\indexB) \: f^{\left( \indexC\right)}(\lambdapoint_\indexB), && \lambdapoint_\indexB \in \left[\uu_\indexB, \uu_{\indexB+\pu+1} \right], \\
	\label{eq:deBoorFixNewtonBasis}
	\Nbasis_{\indexB,\pu}(\uu) &= \prod_{\indexD=1}^\pu \left( \uu - \uu_{\indexB+\indexD}  \right).
\end{align}
Moreover, it can be proven that for all $\indexG$,
\begin{align}
	 \lambda_{\indexG,\pu}( \Bspline_{\indexB,\pu} ) = \delta_{\indexG\indexB}  && \und && \lambda_{\indexG,\pu} \left( \sum_{\indexB=0}^{\totalB-1} \Bspline_{\indexB,\pu} \: \CPara_\indexB   \right) =  \CPara_\indexG\,.
\end{align}

Note that in \myeqref{eq:deBoorFixFunctional}, the evaluation point $\lambdapoint_\indexB$ can be chosen arbitrarily within $\left[\uu_\indexB, \uu_{\indexB+\pu+1} \right]$.
In the following, we will take advantage of this ability and the local nature of this explicit construction of B-spline coefficients $\CPara_\indexA$.

  \section{Extended B-splines}
\label{sec:ExtendedBsplines}

The purpose of taking extended B-splines is to re-establish the stability of a trimmed B-spline basis.
In general, three different types of functions occur if a basis is trimmed: stable, degenerated, and exterior.
They are labeled by $\Bspline_{\indexA,\pu}, \Bspline_{\indexB,\pu},$ and $\Bspline_{\indexC,\pu}$, respectively. 
The latter can be omitted in the following since their support is completely outside of the valid domain, i.e.~$\supp{ \{\Bspline_{\indexC,\pu} \} } \notin \patchdomain$.

\begin{figure}[thb]
  \centering
  \begin{minipage}[b]{0.45\textwidth}
	\tikzsubfig{tikz/conceptExBsplineA.tex}{1.0}{Initial Basis}{fig:extendedBsplineConcept_A}{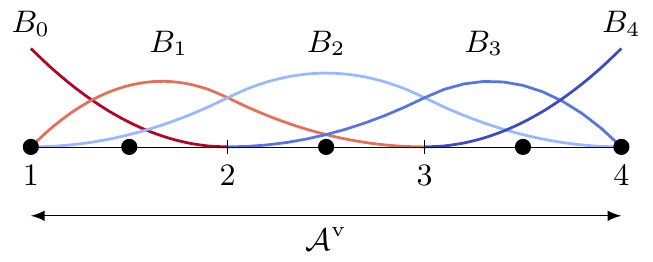}
  \end{minipage}
  \begin{minipage}[b]{0.45\textwidth}
	\tikzsubfig{tikz/conceptExBsplineB.tex}{1.0}{Detection of Degenerated B-splines}{fig:extendedBsplineConcept_B}{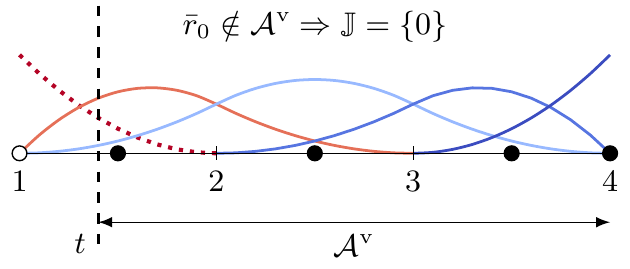}
  \end{minipage}
  \begin{minipage}[b]{0.45\textwidth}
	\tikzsubfig{tikz/conceptExBsplineBB.tex}{1.0}{Extensions of Polynomial Segments}{fig:extendedBsplineConcept_C}{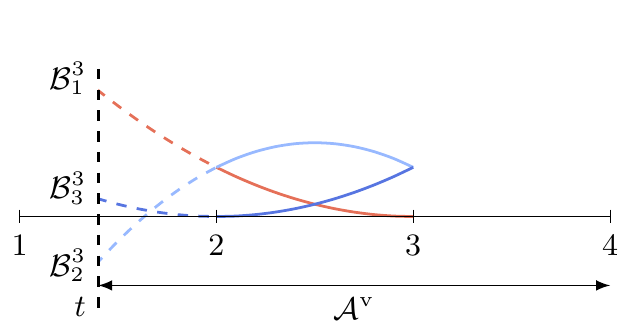}
  \end{minipage}
  \begin{minipage}[b]{0.45\textwidth}
	\tikzsubfig{tikz/conceptExBsplineC.tex}{1.0}{Extended B-splines}{fig:extendedBsplineConcept_D}{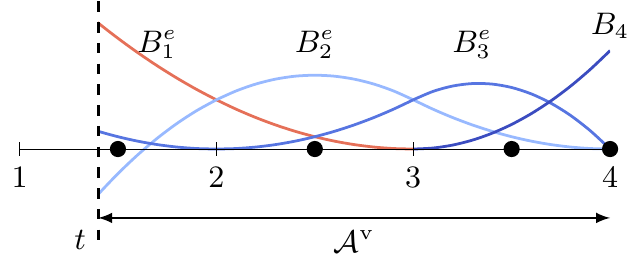}
  \end{minipage}  
  \caption{
  Basic procedure to get from (a) conventional to (d) extended B-splines: (b)~determination of degenerated B-splines and substitution of trimmed polynomial segments by (c) extensions of non-trimmed ones.
  }
  \label{fig:extendedBsplineConcept}
\end{figure}

The basic concept of extended B-splines~$\exBspline_{\indexA,\pu}$ is illustrated in \myfigref{fig:extendedBsplineConcept}.~%
The original basis \ref{fig:extendedBsplineConcept_A} is defined by the knot vector $\KV = \{1,1,1,2,3,4,4,4\}$.
The valid domain of the trimmed patch $\patchdomain$ is determined by a trimming parameter $\trim$, i.e.~$\patchdomain \in (\trim,4]$, as indicated by the dashed vertical line.
In the initial step, the stable B-splines $\Bspline_{\indexA,\pu}$ and degenerated B-splines~$\Bspline_{\indexB,\pu}$ of the basis need to be identified.
Therefore, the size of the supports inside $\patchdomain$ may be evaluated \cite{Hoellig2003b}.
Here, a different approach is proposed: B-splines are labeled degenerated if their support is partially within, but the corresponding anchor is outside of the valid domain, i.e.~$\abscissa_{\indexB} \notin \patchdomain$. The related indices are stored in the index-set~$\degSet$.
In the given example this is the case for $\Bspline_0$.
Next, the polynomial segments of the trimmed knot span are substituted by \emph{extensions} of the polynomial segments $\BsplineSeg^{\indexSpan}_{\indexA}$ of the closest non-trimmed knot span $\indexSpan$ which contains stable B-splines only.  
Finally, those extensions are defined by a linear combination of the degenerated B-splines and the corresponding stable one
\begin{align}
	\label{eq:extendedBspline}
	\exBspline_{\indexA,\pu} = \Bspline_{\indexA,\pu} + \sum_{\indexB \in \degSet_{\indexA}} \eweight_{\indexA,\indexB} \Bspline_{\indexB,\pu} 
\end{align}
where $\degSet_{\indexA}$ is the index-set of all degenerated B-splines related to the current $\exBspline_{\indexA,\pu}$.
Note that \myeqref{eq:extendedBspline} defines $\exBspline_{\indexA,\pu}$ outside of the trimmed knot span, where $\exBspline_{\indexA,\pu} = \Bspline_{\indexA,\pu}$, as well.
It remains to define the extrapolation weights $\eweight_{\indexA,\indexB}$.

\subsection{Univariate Extrapolation Weights}
\label{sec:ExtendedBsplines1Dweights}

In order to compute the extrapolation weights $\eweight_{\indexA,\indexB}$ of an extended B-spline, an interpolation problem needs to be solved.
In particular, the extended polynomial segments $\BsplineSeg^{\indexSpan}_\indexA$ of the non-trimmed knot span~$\indexSpan$ shall be approximated by B-splines of the trimmed knot span~$\indexSpanTrimmed$ such that
\begin{align}
	\label{eq:extendedSplineInterpoalation}
	\BsplineSeg^{\indexSpan}_\indexA \operate{\uu} &= \sum_{\indexB=\indexSpanTrimmed-\pu}^{\indexSpanTrimmed}  { \Bspline_{\indexB,\pu}\operate{\uu} \: \eweight_{\indexA,\indexB} } , & \uu \in \left[ \uu_{\indexSpanTrimmed},\uu_{\indexSpanTrimmed+1} \right).
\end{align}
It should be noted that $\BsplineSeg^{\indexSpan}_\indexA$ can be exactly represented by $\Bspline_{\indexB,\pu}$ since they are polynomials within the \mbox{B-spline} space $\BsplineSpace_{\pu,\KV}$.
Moreover, the coefficient $\eweight_{\indexA,\indexA}$ must be equal to \num{1} due to the fact that $\BsplineSeg^{\indexSpan}_\indexA\operate{\uu} \equiv \Bspline_{\indexA,\pu}\operate{\uu}$, $\uu \in \left[ \uu_{\indexSpan},\uu_{\indexSpan+1} \right)$.
Analogously, 
$\eweight_{\indexA,\indexG}$ must be equal to zero
for other non-zero stable B-splines~$\Bspline_{\indexG,\pu}$ of the knot span~$\indexSpan$, 
since $\BsplineSeg^{\indexSpan}_\indexA\operate{\uu} \neq \Bspline_{\indexG,\pu}\operate{\uu}$, $\uu \in \left[ \uu_{\indexSpan},\uu_{\indexSpan+1} \right)$.
Other values are obtained only if the basis function is degenerated, i.e.~$\indexB \in \degSet_{\indexA}$.
In other words, only certain extrapolation weights have to be computed explicitly.
However, spline interpolation requires anchors $\abscissa$ for each basis function involved as described in~\mysecref{sec:splineInterpolation}.
Unfortunately, the recommended Greville or Demko abscissae are not in general within the trimmed knot span~$\indexSpanTrimmed$. 
Thus, a \emph{quasi interpolation} approach is preferred.
In addition, the proposed scheme introduced in \mysecref{sec:splineInterpolation} is applied to single intervals and hence it is ideally suited to compute the extrapolation weights $\eweight_{{\indexA},{\indexB}}$.

In order to obtain $\eweight_{{\indexA},{\indexB}}$ the function $f$ of the dual functional~\myRef{eq:deBoorFixFunctional} is substituted by $\BsplineSeg^{\indexSpan}_{\indexA}$, leading to
\begin{align}
	\label{eq:EW}
	\eweight_{{\indexA},{\indexB}} &= \lambda_{\indexB,\pu}(\BsplineSeg^{\indexSpan}_{\indexA}) =
	\frac{1}{\pu!} \sum^\pu_{\indexC=0} (-1)^\indexC \Nbasis^{\left(\pu-\indexC\right)}_{\indexB,\pu}(\lambdapoint_\indexB) \: \BsplineSeg^{\indexSpan^{\left( \indexC\right)}}_{\indexA}(\lambdapoint_\indexB), && \lambdapoint_\indexB \in \left[\uu_\indexB, \uu_{\indexB+\pu+1} \right] .
\end{align}

In fact, \myeqref{eq:EW} provides the entire information required for setting up the extrapolation weights.
However, it may seem a bit complex at first glance.
With this issue in mind, the
evaluation is discussed in more detail.
The polynomial $\Nbasis_{\indexB,\pu}$ can be 
written explicitly in power basis form
\begin{align}
	\label{eq:Nbasis_explicit}
	\Nbasis_{\indexB,\pu}(\uu) &=  \sum^{\pu}_{\indexC=0} \beta_\indexC \: \uu^{\indexC}.
\end{align}
The corresponding coefficients $\beta_\indexC$ are computed by
\begin{align}
	\label{eq:Nbasis_explicit_beta}
	 \beta_\indexC & =   (-1)^\indexC
	 \sum_{\indexG=1}^\totalG  \prod_{\indexD \in \NSet_{\indexC,\indexG}} \uu_\indexD 
	 && \with &&
	 \totalG = \frac{ \pu! }{ \left( \pu - \indexC \right)! \: \indexC! }
\end{align}
where the sum over $\NSet_{\indexC,\indexG}$ denotes all $\indexC$-combinations with repetition
of the knots appearing in the definition~\myRef{eq:deBoorFixNewtonBasis} of $\Nbasis_{\indexB,\pu}$, i.e.~$ \uu_{\indexB+1},\dots,\uu_{\indexB + \pu}$.
The segments $\BsplineSeg^{\indexSpan}_{\indexA}$ of $\Bspline_{\indexA,\pu}$ can be expressed by a Taylor expansion
\begin{align}
	\label{eq:Bspline_Taylor}
	\BsplineSeg^{\indexSpan}_{\indexA}(\uu)
	& = \sum^{\pu}_{\indexC=0} \frac{\Bspline_{\indexA,\pu}^{(\indexC)}(\taylorpoint)}{\indexC!}\:(\uu-\taylorpoint)^{\indexC}
	= \sum^{\pu}_{\indexC=0} \alpha_\indexC\:(\uu-\taylorpoint)^{\indexC} ,  
	&& \taylorpoint \in \left[\uu_{\indexSpan}, \uu_{\indexSpan+1}\right)
\end{align}
where the point~$\taylorpoint$ is within the corresponding knot span $s$.
\myEqref{eq:Bspline_Taylor} can also be written in power basis form such that
\begin{align}
	\label{eq:Bspline_PowerBasis}
	\BsplineSeg^{\indexSpan}_{\indexA}(\uu) &=
	\sum^{\pu}_{\indexC=0} \poly_\indexC \: \uu^\indexC 
	&& \with && 
	\poly_\indexC = \sum_{\indexD=\indexC}^{\pu} \binom{\indexD}{\indexC} \alpha_\indexD \left(-\taylorpoint \right)^{\indexD-\indexC}
\end{align}
where $\binom{\indexD}{\indexC}$ denotes the binomial coefficient defined as
\begin{align}
	\label{eq:binomialcoefficient}
	\binom{\indexD}{\indexC} \coloneqq \frac{\indexD!}{(\indexD-\indexC)! \:\indexC!} \: .
\end{align}

The resulting extrapolation weight $\eweight_{\indexA,\indexB}$ does not change if the location of the evaluation point $\lambdapoint_\indexB$ is modified.
Thus, functional values depending on $\lambdapoint_\indexB$ cancel out \cite{Rueberg2014notes}.
Consequently, the dual functional can be written as
\begin{align}
	\label{eq:deBoorFix_explicit}
		\eweight_{{\indexA},{\indexB}} &= \lambda_{\indexB,\pu}(\BsplineSeg^{\indexSpan}_{\indexA}) =
	\frac{1}{\pu!} \sum^\pu_{\indexC=0} (-1)^\indexC 
	\left( \pu - \indexC \right)! \: \beta_{\pu-\indexC} \:
	 \indexC ! \:  \poly_\indexC
\end{align}
where $\eweight_{{\indexA},{\indexB}}$ depends only on the coefficients $\beta$ and $\poly$ of explicit representations \myRef{eq:Nbasis_explicit} and \myRef{eq:Bspline_PowerBasis} of the polynomials.
Supplementary information regarding the evaluation of $\lambda_{\indexB,\pu}(\BsplineSeg^{\indexSpan}_{\indexA})$ is provided in~\myappref{chap:extrapolationWeightsExample}.
In case of a uniform knot vector, a simplified formula can be derived  which solely relies on the indices of the  B-splines involved, see e.g.~\cite{Hoellig2003b}.

\subsection{Bivariate Extrapolation Weights}

\begin{figure}[htb]
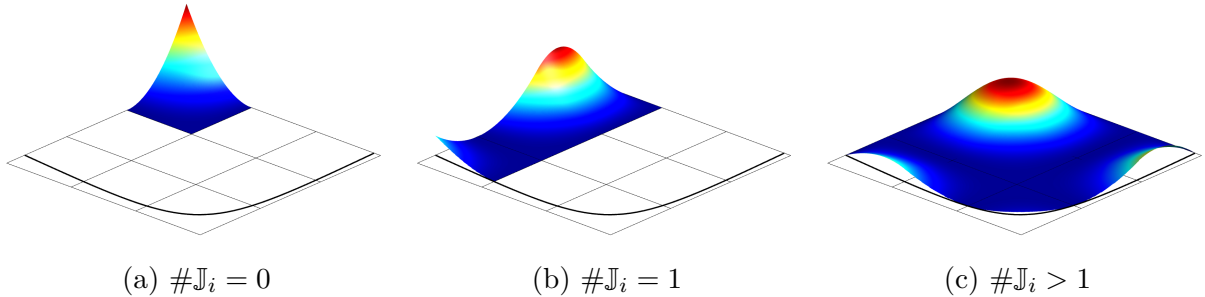

  \centering
  \begin{minipage}[b]{0.32\textwidth}
	\tikzsubfig{tikz/exBasisPlot2.tex}{0.13}{$\myCardinality{\degSet_{\indexA}} = 0$}
	{fig:extendedBsplinePlot2}{Figure5a}
  \end{minipage}
  \begin{minipage}[b]{0.32\textwidth}
	\tikzsubfig{tikz/exBasisPlot3.tex}{0.13}{$\myCardinality{\degSet_{\indexA}} = 1$}
	{fig:extendedBsplinePlot3}{Figure5b}
  \end{minipage}
  \begin{minipage}[b]{0.32\textwidth}
	\tikzsubfig{tikz/exBasisPlot1.tex}{0.13}{$\myCardinality{\degSet_{\indexA}} > 1$}
	{fig:extendedBsplinePlot1}{Figure5c}
  \end{minipage}
  \caption{Bivariate extended B-splines $\exBspline_{\indexA,\pu}$ with various 
  cardinalities of the index-set $\degSet_{\indexA}$ which indicates the number of related degenerated B-splines. Note that (a) is in fact a conventional B-spline, i.e.~$\exBspline_{\indexA,\pu} \equiv \Bspline_{\indexA,\pu}$, since $\degSet_{\indexA}$ is empty.}
  \label{fig:extendedBsplineBasisPlots}
\end{figure}

Bivariate extrapolation weights $\eweight_{{\indexA},{\indexB}}$ are obtained by the tensor product of their univariate counterparts calculated for each parametric direction.
The construction procedure is visualized in~\myfigref{fig:BivariateExtrapolationWeights} and examples of different bivariate extended B-splines are shown in~\myfigref{fig:extendedBsplineBasisPlots}.%
\tikzfigposition{tikz/Eij_2d.tex}{1.6}{The construction of bivariate extrapolation weights $\eweight_{{\indexA},{\indexB}}$. The basis is defined by a tensor product of the B-spline depicted in~\myfigref{fig:extendedBsplineConcept_A}, the dashed line indicates the trimming curve and $\patchdomain$ is highlighted in gray. Stable B-splines are marked by black and green circles. The shown values of $\eweight_{{\indexA},{\indexB}}$ are related to the degenerated basis function marked by the blue circle in the upper right corner of the parameter space. The B-splines of the closest non-trimmed knot span are indicated by green circles. The values of the univariate extrapolation weights are derived in \myappref{chap:extrapolationWeightsExample}.}{fig:BivariateExtrapolationWeights}{Figure4}{t}
The closest non-trimmed knot span can be detected by measuring the distance between the anchor of the degenerated B-spline and the midpoints of the surrounding knot spans.
In general, a degenerated B-spline is distributed to $\left(\pusurf+1\right)\left(\pvsurf+1\right)$ stable ones.
Further, several degenerated B-splines $\Bspline_{\indexB,\pu}$ may be associated to a single stable B-spline~$\Bspline_{\indexA,\pu}$.
In particular, the number of $\Bspline_{\indexB,\pu}$ related to $\Bspline_{\indexA,\pu}$ is determined by the cardinality of the corresponding index-set $\myCardinality{\degSet_{\indexA}}$.
However, the extension procedure is restricted to those basis functions which are close to the trimming curve.
For instance, the basis function shown in \myfigref{fig:extendedBsplinePlot2} is a conventional B-spline since it is far enough away from the trimming curve.
The actual size of the affected region depends on the fineness of the parameter space and the degree of its basis functions. 

\newpage
\subsection{Properties of Extended B-splines}

Extended B-splines inherit most essential properties of conventional B-splines~\cite{Hoellig2003b,Hoellig2003a}.
First of all, they are \emph{linear independent} and \emph{polynomial precision} is guaranteed. 
Thus, they form a basis for a spline space.
Each knot span has exactly $\pu + 1$ non-vanishing basis functions which span the space of all polynomial of degree~$\leqslant \pu$ over $\patchdomain$. 
Extended B-splines have \emph{local support} since only B-splines near the trimming curve are affected by the extension procedure.
In addition, approximation estimates have the \emph{same convergence order} as conventional B-splines.

However, there are also some differences.
It is important to note that the extrapolation weights may be \emph{negative}, 
hence the evaluation of extended B-splines may lead to negative values.
Conventional B-splines, on the other hand, are strictly non-negative.
This property is exploited in some contact formulations~\cite{Temizer2011a} and structural optimization~\cite{Nagy2010a}, for instance. 
In such cases, the application of extended B-splines requires further considerations.
The main difference in favor of extended B-splines is the \emph{stability} of the corresponding basis.
In particular, the condition number of a system is independent of the location of the trimming curve due to the substitution of B-splines with small support. 
This is an important advantage since it guarantees the performance of iterative solvers and robust solutions. Theoretically, this also ensures larger stable step sizes in explicit time integration schemes leading to faster simulations.
Another benefit of extended B-splines is that the determination of proper \emph{anchors} for trimmed parameter spaces is straightforward since all anchors are within $\patchdomain$ by construction.
This is an essential feature for interpolation and collocation schemes.

  \section{Analysis of Trimmed Geometries}
\label{sec:analysis}

The extended B-spline basis is applied to an isogeometric BEM framework.
Firstly, the main points of an isogeometric boundary element formulation are briefly introduced.
Then, the partitioning of integration elements over trimmed surfaces is outlined. 
Finally, the set up of the \emph{stable} system of equations is discussed.

It is noteworthy that the stabilization of a trimmed basis by extended B-splines also applies for finite element simulations.
At this point, the distinguishing feature of the BEM is that \emph{all} boundary conditions are employed in an integral sense.
In case of finite elements, the enforcement of essential boundary conditions requires additional attention, see e.g.~\cite{Hoellig2003b,Rueberg2012a}.

\subsection{Isogeometric Boundary Element Analysis}
\label{sec:bem}

The direct \emph{boundary integral equation}
\begin{align}
    \label{eq:bie}
     c\ofpt{x} \primary\ofpt{x} &= \int_\boundary \fund{U}(\pt{x},\pt{y}) \: \dual \ofpt{y} \dgamma{y}
    - \int_\boundary\fund{T}(\pt{x},\pt{y}) \:  \primary\ofpt{y} \dgamma{y} 
    & &\forall \pt{x}, \pt{y}\in \boundary
\end{align}
provides the basis of the boundary element method with $\boundary$ representing the boundary of a body $\domain$ subjected to external loading without body forces.
The jump term $c\ofpt{x}$ depends on the geometrical angle of $\boundary$ at $\pt{x}$ and is $0.5$ where $\boundary$ is smooth.
The primary and dual field variables are denoted by $\primary$ and $\dual$, respectively.
On the other hand, the terms $\fund{U}$ and $\fund{T}$ refer to the corresponding fundamental solutions.
In general, a fundamental solution $\fund{U}(\pt{x},\pt{y})$ provides the response at a \emph{field point} $\pt{y}$ due to a unit point source applied at $\pt{x}$, which is denoted as \emph{source point}. 
The fundamental solutions for potential and elasticity problems considered in this paper are well known and are given in various textbooks, see e.g.~\cite{BeerSmithDuenser2008b,gaul2003,Sauter2011b,Steinbach2008b}.
The fundamental solutions become singular as $\pt{y}$ approaches $\pt{x}$ and their evaluation requires appropriate regularization and integration schemes as described in \cite{beer2015bb,Marussig2014aCMAME}.

Once \myeqref{eq:bie} is solved, the Cauchy data $\primary\ofpt{y}$ and $\dual\ofpt{y}$ are known and $\primary\ofpt{x}$ in the domain can be obtained by the representation formula 
\begin{align}
  \label{eq:bie:representation-formula}
  \primary\ofpt{x} &= \int_\Gamma \fund{U}(\pt{x},\pt{y}) \: \dual  \ofpt{y} \dgamma{y}
  - \int_\Gamma \fund{T}(\pt{x},\pt{y}) \: \primary\ofpt{y} \dgamma{y} 
    & &\forall \pt{x}\in \domain,~\forall \pt{y}\in \Gamma.
\end{align}

Following the isogeometric concept, the geometry as well as the Cauchy data are discretized by B-splines.
Hence, the corresponding coefficients are expressed by control variables rather than nodal values on $\boundary$.
The \emph{unknown} values are determined by solving a system of equations which is set 
up by enforcing the boundary integral equation at a specific set of \emph{collocation points} $\collocationPoint$.
Their position is determined by the anchors of the related basis functions, i.e.~$\collocationPoint_\indexA = \mychi\operate{\greville_\indexA}$.
The resulting discretized boundary integral equation for the Neumann problem is given by
\begin{align}
    \label{eq:BIEdisNeumann}
    \myMat{K} \tilde{\myVec{\primary}} & = 
	\myMat{V} \tilde{\myVec{g}}_N.
\end{align}
The vector $\tilde{\myVec{\primary}}$ contains the unknown control variables of the Dirichlet field while the known Neumann coefficients are collected in $\tilde{\myVec{g}}_N$.
The matrix $\myMat{K}$ is singular for an interior Neumann problem, but invertible in case of an exterior problem \cite{Sauter2011b}.
The entries of the right hand side matrix $\myMat{V}$ and the left hand side matrix $\myMat{K}$ are determined by
\begin{align}
  \label{eq:BIEdisEntries}
  \myMat{V}\left[\indexA,\indexB\right] & = 
  \int_{\boundary} \fund{U}(\collocationPoint_\indexA,\pt{y}) \DualBasis_\indexB\ofpt{y} \dgamma{y}
  && \und && \myMat{K}\left[\indexA,\indexB\right] = 
  {c}_{ij} + 
  \int_{\boundary} \fund{T}(\collocationPoint_\indexA,\pt{y}) \PrimaryBasis_\indexB\ofpt{y} \dgamma{y}
\end{align}
where $\DualBasis_\indexB$ and $\PrimaryBasis_\indexB$ are the univariate or bivariate B-spline basis function used for the discretization of $\dual$ and~$\primary$. 
The jump term coefficients ${c}_{\indexA\indexB}$ are calculated such that 
\begin{align}
    \sum_{\indexB=0}^{\totalB} {c}_{\indexA\indexB} \: \PrimaryBasis_\indexB\operate{\collocationPoint_\indexA} =  c \operate{\collocationPoint_\indexA}.
\end{align}
Hence, ${c}_{\indexA\indexB}$ is non-zero only if $\collocationPoint_\indexA$ is in the support of the basis function $\PrimaryBasis_\indexB$.
In case of a mixed problem the equations are sorted with respect to the Dirichlet and Neumann boundary which leads to a block system of matrices~\cite{Marussig2014aCMAME}. 
In order to perform the integral over $\boundary$, an appropriate partition of integration elements is required, which will be discussed for trimmed geometries in the following section.

\subsection{Integration Elements}

For non-trimmed knot spans the definition of integration elements is straightforward, but trimmed ones require a different treatment. 
The detection of trimmed knot spans can be performed accordingly to \citet{Schmidt2012a}:
\begin{enumerate}
	\item Determine \emph{intersection points} of the trimming curve and the grid produced by the tensor product of the knot vectors.
	\item Detect invalid \emph{cutting patterns} and perform knot insertion to obtain valid ones.
	\item Assign an \emph{element type} to each knot span based on the valid cutting patterns.
\end{enumerate}

\begin{figure}[htb]
  \centering
  \begin{minipage}[b]{0.24\textwidth}
	\tikzsubfig{tikz/trimmingCuttingPatternA.tex}{2.0}{}{fig:CuttingPatternTriangle}{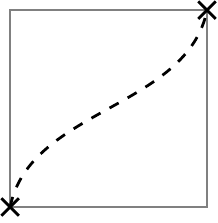}
  \end{minipage}
  \begin{minipage}[b]{0.24\textwidth}
	\tikzsubfig{tikz/trimmingCuttingPatternB.tex}{2.0}{}{fig:CuttingPatternTriangleQuad}{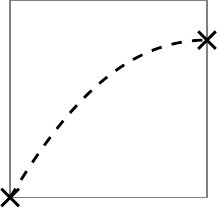}
  \end{minipage}
  \begin{minipage}[b]{0.24\textwidth}
	\tikzsubfig{tikz/trimmingCuttingPatternC.tex}{2.0}{}{fig:CuttingPatternQuad}{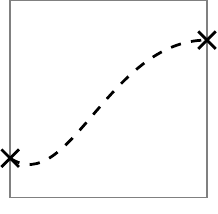}
  \end{minipage}
  \begin{minipage}[b]{0.24\textwidth}
	\tikzsubfig{tikz/trimmingCuttingPatternD.tex}{2.0}{}{fig:CuttingPatternPentagon}{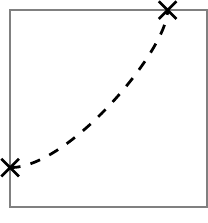}
  \end{minipage}
  \caption{Illustration of valid cutting patterns of a single knot span. The actual element type is determined by the direction of the trimming curve which is indicated by the dashed line. The intersection points are highlighted by crosses.}
  \label{fig:trimmingCuttingPattern}
\end{figure}

\myFigref{fig:trimmingIntersectionpoints} shows a trimmed parameter space with related element types and \myfigref{fig:trimmingCuttingPattern} depicts the valid cutting patterns considered in this paper. 
Elements of type $1$ are the regular knot spans, whereas elements of type $-1$ are not considered during the analysis since they are completely outside of the valid domain.
The shape of trimmed knot spans corresponds to either a triangle, a quadrilateral or a pentagon.%
\tikzfigposition{tikz/trimmingIntersectionpoints.tex}{1.6}{
Trimmed parameter space and corresponding element types: \num{1} labels untrimmed knot spans whereas \num{-1} denotes knot spans which are outside of the computational domain. In case of a trimmed knot span the element type indicates the number of its edges in $\patchdomain$, i.e.~\num{3} (triangle), \num{4} (quadrilateral) or \num{5} (pentagon).
}{fig:trimmingIntersectionpoints}{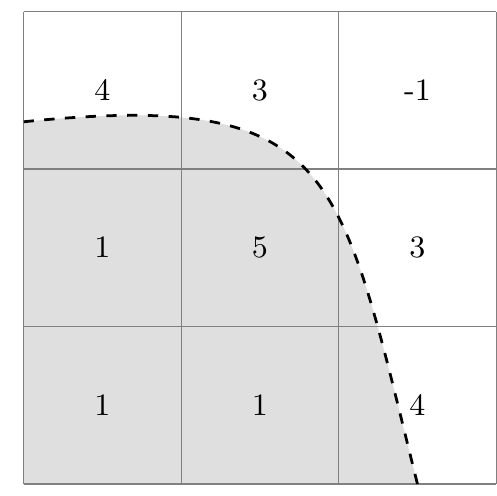}{hb}%
Accordingly, they are labeled as elements of type $3$, $4$, and $5$.
It remains to define the mapping~$\mychi_{\xi}\ofpt{\xi}$ from the reference element~$\ieReference=\left[-1,1\right]^{d-1}$ to elements of type $\{3,4,5\}$.

\begin{figure}[thb]
  \centering
  \begin{minipage}[b]{0.47\textwidth}
	\tikzsubfig{tikz/integrationRegionsRegular.tex}{1.7}{Regular Local Coons Patch}
	{fig:localCoonsRegular}{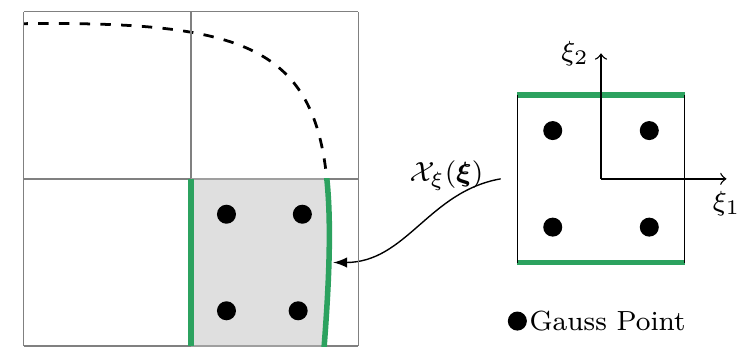}
  \end{minipage}
  \hfill
    \begin{minipage}[b]{0.47\textwidth}
	\tikzsubfig{tikz/integrationRegionsDegenerated.tex}{1.7}{Degenerated Local Coons Patch} 
	{fig:localCoonsDegernated}{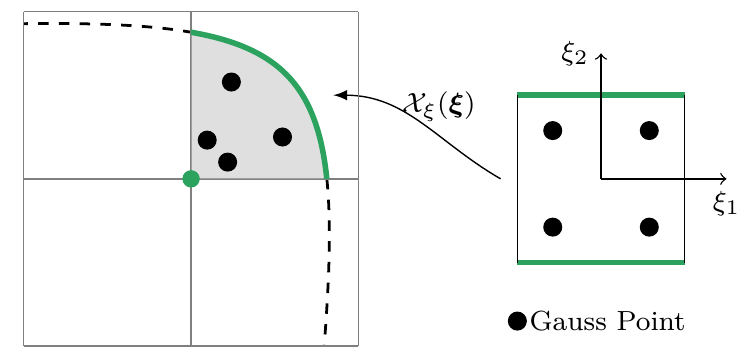}
  \end{minipage}
    \caption{Mapping $\mychi_{\xi}\ofpt{\xi}$ from the reference element to the B-spline parameter space. The related knot span is marked in gray, the dashed line represents the trimming curve, circles indicate quadrature points and higher degree edges, i.e. $\pu > 1$, are highlighted in green.}
  \label{fig:localCoonsPatches}
\end{figure}

There are various approaches to define the local mapping $\mychi_{\xi}\ofpt{\xi}$, e.g.~\cite{Fries2015a,Kim2009a,Schmidt2012a,Wang2013a}.
In this paper, the idea of \cite{Beer2015a} is locally applied to the trimmed knot spans.
By adapting the concept of Coons patches, the edges of the trimmed knot span are used to specify $\mychi_{\xi}\ofpt{\xi}$. 
There is a single edge determined by the trimming curve, the others are straight lines related to the grid of the parameter space.
The latter can be represented by the grid points and linear B-splines $\Bspline_{\indexB,1}$ based on the knot vector $\KV=\left\{-1,-1,1,1\right\}$.
The former is obtained by the application of knot insertion so that the trimming curve becomes interpolatory at the intersection points.
The resulting segment $\curveBez$ and the opposing straight edge~$\curveBezOpposite$ determine $\mychi_{\xi}\ofpt{\xi}$ by the following construction:
firstly, it is ensured that both curves are defined on the same parameter range, i.e.~$\xi_1 \in \left[-1,1\right]$.
Next, the knot vector~$\KV_e$ of $\curveBezOpposite$ is refined by degree elevation and knot insertion so that it is equivalent to the knot vector~$\KV_b$ of $\curveBez$.
Consequently, both edges are described by the same basis functions $\Bspline_{\indexA,\pu}\operate{\xi_1}$.
Combined with a linear interpolation $\Bspline_{\indexB,1}\operate{\xi_2}$ given by the knot vector $\KV=\left\{-1,-1,1,1\right\}$ the integration region of the trimmed knot span is represented by
\begin{align}
	\label{eq:ruledSurface_mapping}
	\mychi_{\xi}\ofpt{\xi} &\coloneqq \pt{x}\operate{\xi_1,\xi_2} = 
	\sum_{\indexA=0}^{\totalA-1} \sum_{\indexB=0}^{1} \Bspline_{\indexA,\pusurf} \operate{\xi_1}  \Bspline_{\indexB,1}\operate{\xi_2}    \: \CPCoons_{\indexA,\indexB}
\end{align}
where $\CPCoons_{\indexA,\indexB}$ denote the control points of $\curveBez$ and $\curveBezOpposite$ within the parametric space.
In case of the element type \num{3} the opposite edge $\curveBezOpposite$ degenerates to a point.
Moreover, an element of type \num{5} can be treated by subdivision into three elements of type 3.
The local mapping $\mychi_{\xi}\ofpt{\xi}$ is exemplary shown in \myfigref{fig:localCoonsPatches}.
It is noteworthy that the integration points do not coincide at the degenerated point since an open Gauss quadrature rule is used.

\newpage
\subsection{Stabilization}

The approach presented in the previous section permits a proper integration over the valid area $\patchdomain$ of a patch, since the partition of integration elements $\ieLocal$ is such that
\begin{align}
  \patchdomain = \bigcup_{\indexA=1}^{\totalA} \ieLocal_\indexA .
\end{align}
Moreover, it allows the integration over the boundary $\boundary$ which is defined by the union of all $\patchdomain$. 
Basis functions need to be integrated within each $\ieLocal$
in order to set up the system of equations for the analysis. 
If extended B-splines $\exBspline_{\indexA,\pu}$ are involved, they could be evaluated directly by \myeqref{eq:extendedBspline}.
However, they can be included into the simulation without being considered during the integration process as well.
In particular, integration is performed with the original basis functions~$\Bspline_{\indexA,\pu}$,  $\indexA \notin \exteriorSet$, where the index-set $\exteriorSet$ denotes all exterior B-splines~$\Bspline_{\indexC,\pu}$ which are completely outside of the domain, i.e.~$\supp \{\Bspline_{\indexC,\pu} \}  \notin \patchdomain$.
In other words, the computation of entries in the system matrix does not differ from the regular case of non-trimmed patches despite the tailored mapping $\mychi_{\xi}\ofpt{\xi}$ from the reference element to the parameter space which is applied for trimmed knot spans.

Consequently, the resulting system matrix $\myMat{K} \in \R^{ n \times m }$ is rectangular, where $n$ denotes the number of collocation points and $m$ represents the total number of integrated $\Bspline_{\indexA,\pu}$,  $\indexA \notin \exteriorSet$.
In order to obtain a square matrix an \emph{extension matrix} $\myMat{E} \in \R^{ m \times n }$ is introduced \cite{Hoellig2003b}.
This sparse matrix $\myMat{E}$ contains all extrapolation weights $\eweight_{{\indexA},{\indexB}}$ including the trivial ones, i.e.~$\eweight_{{\indexA},{\indexA}} = 1$.
The transformation of the original to the stable extended B-spline basis is performed by multiplying the extension matrix to the left hand side and ride hand side of the system of equations. 
In context of an exterior Neumann problem the stabilization is given by
\begin{align}	
      \label{eq:extensionOperatorApplication}
	\myMat{K} \myMat{E}_\PrimaryBasis \tilde{\myVec{\primary}} & = 
	\myMat{V} \myMat{E}_\DualBasis \tilde{\myVec{g}}_N    \\
	\label{eq:stableSystem}
	 {\myMat{K}_e}  \tilde{\myVec{\primary}} & = \myVec{f} 
\end{align}
where $\tilde{\myVec{g}}_N \in \R^{n}, \tilde{\myVec{\primary}} \in \R^{n}$, $\myVec{f} \in \R^{n}$ and ${\myMat{K}_e}\in \R^{ n \times n }$. 
The subscripts of $\myMat{E}$ emphasize that the extrapolation weights are related to the basis functions $\PrimaryBasis$ and $\DualBasis$ of the Cauchy data. 
The resulting stable system~\myRef{eq:stableSystem} is subsequently solved and the obtained solution $\tilde{\myVec{\primary}}$ corresponds to the extended B-splines of the unknown field.
In case of multi-patch geometries, the extrapolation weights~$\eweight_{{\indexA},{\indexB}}$ of each patch have to be assembled to $\myMat{E}$ with respect to the global degrees of freedom.
However, the application of the extension operator is particularly convenient, if extended B-splines are added to an existing isogeometric code.

  \section{Numerical Results}
\label{sec:results}

The proposed stabilization of trimmed parameter spaces is investigated in this section. 
Firstly, the approximation quality of extended B-splines is examined by means of interpolation problems.
Secondly, the interplay of knot spacing and stabilization is studied.
Afterwards, extended B-splines are applied to isogeometric BEM examples defined by simple trimming cases.
Finally, the applicability to a complex geometry is shown.
In the following, the anchors of the basis functions are generally determined by the Greville abscissae~\myRef{eq:greville}.

\subsection{Assessment of Approximation Quality}

In order to verify the approximation quality of extended B-splines the following test case is studied:
the initial B-spline basis defined by $\KV = \left\{-1,-1,1,1\right\}$ is refined by
degree elevation up to degrees $\pu=\{2,3,4\}$ and uniform knot insertion.
The knot insertion depth~$\knotinsertion$ indicates how often the knot spans have been subdivided.
Hence, it defines how many knots are inserted. 
A trimming parameter~$\trim <1$ determines the valid domain of the patch $\patchdomain \in [-1,\trim)$.

\begin{figure}[b!]
  \centering
  \begin{minipage}[b]{0.45\textwidth}
	\tikzsubfig{tikz/trimmedbsplineExample1d_extendedBspline.tex}{1.7}{Extended B-splines}{fig:exampleSchemes_1dEx}{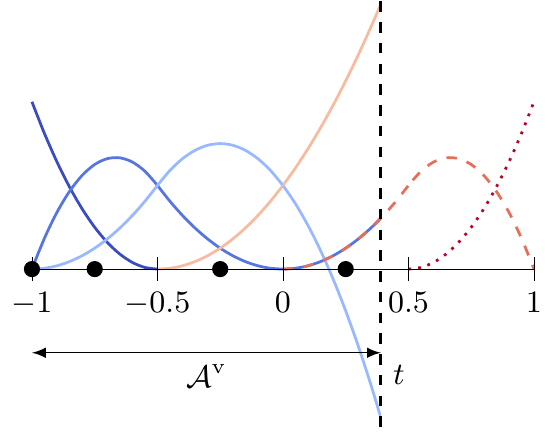}
  \end{minipage}
  \begin{minipage}[b]{0.45\textwidth}
	\tikzsubfig{tikz/trimmedbsplineExample1d_naive.tex}{1.7}{Conventional B-splines}{fig:exampleSchemes_1dNaive}{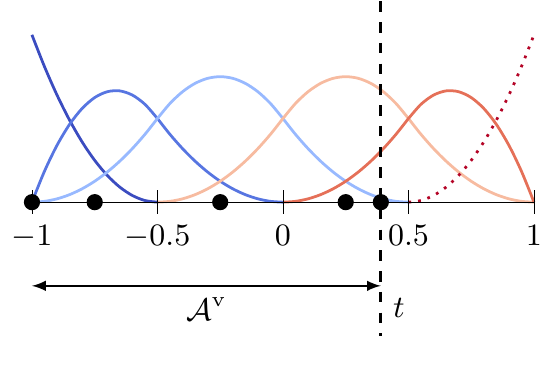}
  \end{minipage}
  \caption{Trimmed univariate basis specified by $\pu=2, \knotinsertion = 2$ and $\trim=0.4$ for (a) extended B-splines and (b) the simple approach using conventional B-splines. The involved basis functions are indicated by continuous lines and their anchors are marked by circles. 
  Dotted basis functions do not contribute to the interpolation problem. In (a), the dashed basis function indicates a degenerated B-spline.
  }
  \label{fig:exampleSchemes}
\end{figure}

An interpolation problem is solved using extended B-splines for each trimmed basis. 
In addition a simple approach is applied, where conventional B-splines are used, but the anchors of degenerated basis functions are merely shifted into $\patchdomain$.
In both cases, a support does not contribute to the approximation if it is entirely outside of the domain.
The schemes are depicted in \myfigref{fig:exampleSchemes}. 
The quality and stability of the approximation $f_h(\uu)$ are specified by the relative interpolation error measured in the $\Ltwo$-norm $\myNorm{ \err_{rel}}_{\Ltwo}$ as well as the condition number of the spline collocation matrix~$\kappa(\myMat{A}_\uu)$.

\subsubsection{Spline Interpolation in 1D}

The target function of the interpolation problem is given by
\begin{align}
	\label{eq:originalFunction1d}
	f(\uu) &= \frac{1}{|a-\uu|} && \with && a = -1.1 \notin \patchdomain,~ 
        \uu \in \patchdomain .
\end{align}	
The knot insertion depth~$\knotinsertion$ is set to \num{4} and the interpolation is performed for several trimming parameters $\trim$, $0.5 < \trim < 1$.
The resulting $\kappa(\myMat{A}_\uu)$ and $\myNorm{ \err_{rel}}_{\Ltwo}$ of the untrimmed case $\trim=1$ are given in~\mytabref{tab:untrimmedResults} for reference purpose. 
All other results are summarized in \myfigref{fig:example1d_cond} and \myfigref{fig:example1d_err}.

\sisetup { round-mode  = places, round-precision = 3 }
\begin{mytableposition}
  {Results of the interpolation problem of the untrimmed basis, i.e.~$\trim=1$.}%
  {untrimmedResults}
  {ccc}{ht}
  \mytableheader{
	$\pu$ & $\kappa$  & $\myNorm{ \err_{rel}}_{\Ltwo}$ }  
	$2$   & \num{2.5}     & \num{1.98946e-02} \\ 
	$3$   & \num{4.30981} & \num{5.73360e-03} \\ 
	$4$   & \num{7.93821} & \num{1.75000e-03} \\ 
\end{mytableposition}%
\tikzfig{tikz/resultsInterpolation1d_cond}{0.95}
{Condition number $\kappa(\myMat{A}_\uu)$ for several degrees $\pu$ and trimming parameters~$\trim$. The labels of the horizontal axis indicate knots of the trimmed basis.}{fig:example1d_cond}{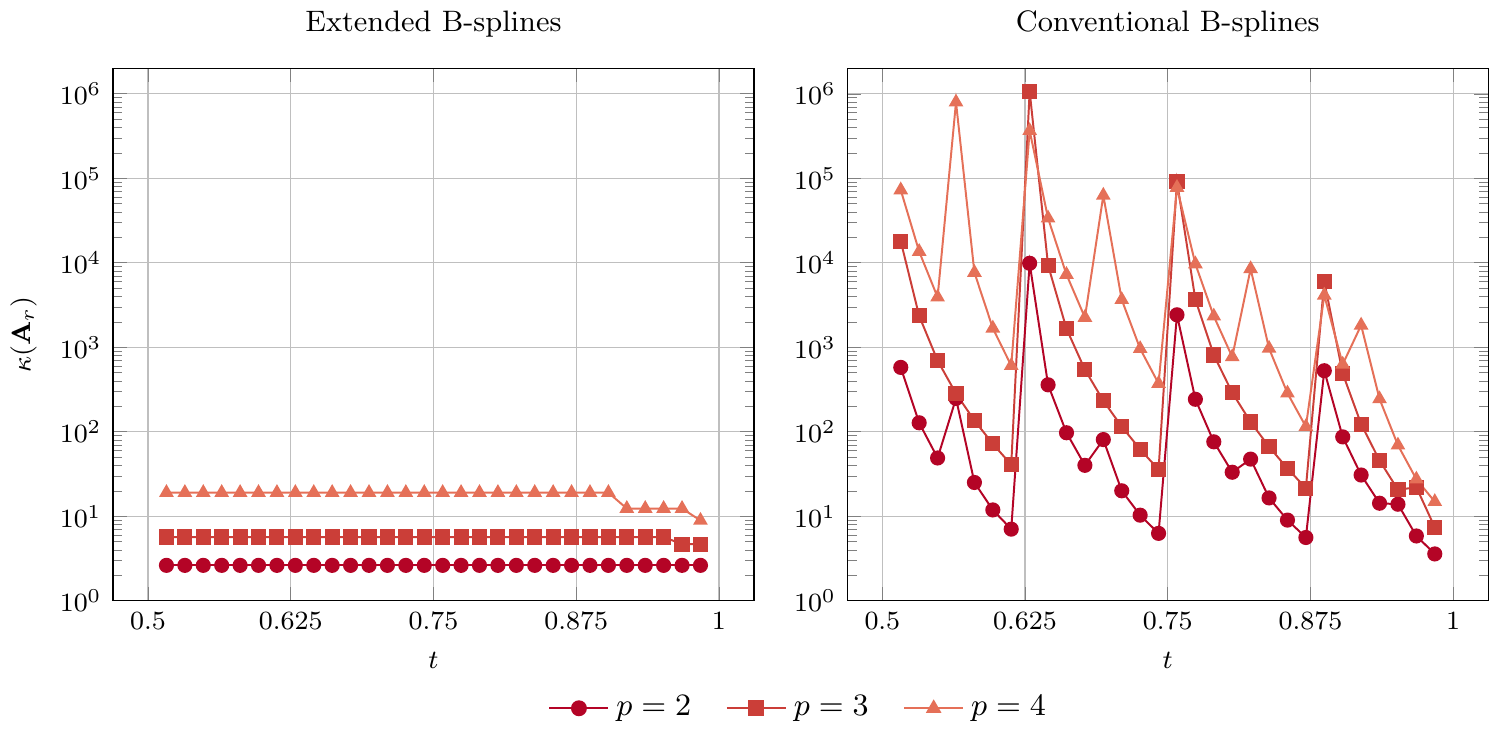}
\tikzfig{tikz/resultsInterpolation1d_err}{0.8}
{Relative interpolation error of the univariate example with several degrees~$\pu$ related to the trimming parameter~$\trim$ (left) and the degrees of freedom $\DOF$ (right). }{fig:example1d_err}{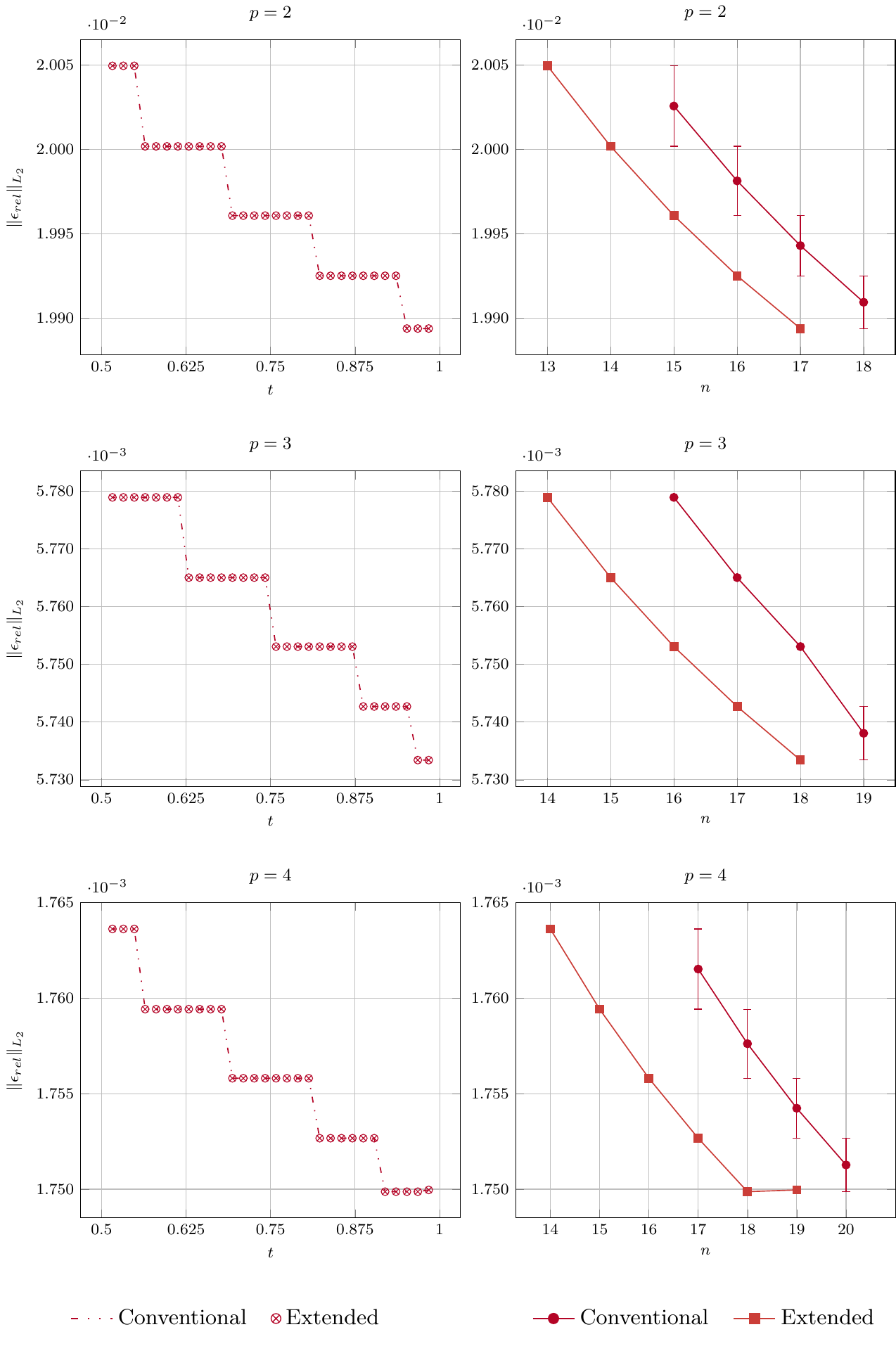}
%

It can be observed that the condition number of $\myMat{A}_\uu$ due to shifted anchors is considerably influenced by the position of $\trim$. 
A peak is reached as soon as $\trim$ approaches an anchor or a knot value. 
If extended B-splines are used, $\kappa(\myMat{A}_\uu)$ hardly changes. 
In other words, the extended B-spline basis is stable. 
At the same time, the relative error is almost identical with respect to $\trim$. 
Moreover, the extended B-spline approach yields even more accurate results, if $\myNorm{ \err_{rel}}_{\Ltwo}$ is related to the number of anchors, i.e.~degrees of freedom $\DOF$. 

\subsubsection{Spline Interpolation in 2D}

\begin{figure}[htb]
  \centering
  \begin{minipage}[b]{0.45\textwidth}
	\tikzsubfigscale{tikz/trimmedbsplineExample2d_extendedBspline.tex}{0.91}{Extended B-splines}{fig:exampleSchemes2d_Ex}{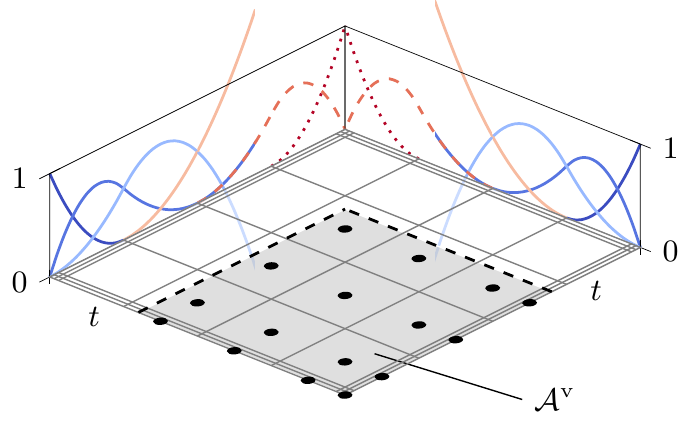}
  \end{minipage}
  \hfill
  \begin{minipage}[b]{0.45\textwidth}
	\tikzsubfigscale{tikz/trimmedbsplineExample2d_naive.tex}{0.91}{Conventional B-splines}{fig:exampleSchemes2d_Naive}{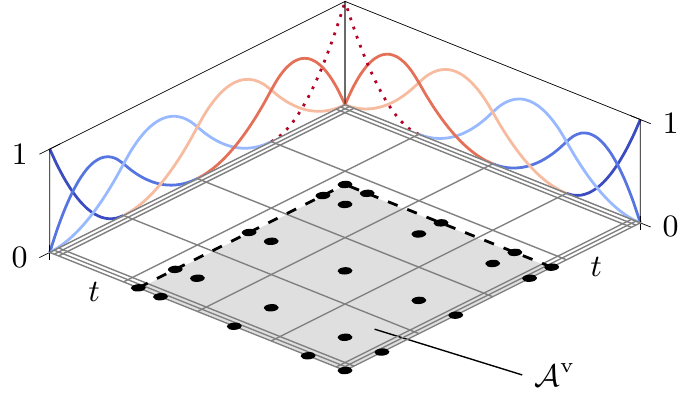}
  \end{minipage}
  \caption{
  Trimmed bivariate basis specified by $\pu=2, \knotinsertion = 2$ and $\trim=0.4$ for (a) extended B-splines and (b) the simple approach using conventional B-splines. 
  The domain~$\patchdomain$ of the interpolation problem is highlighted in gray and circles indicate the basis functions involved.
  Note that there are anchors on the trimming curve for the simple approach.
  }
  \label{fig:exampleSchemes2d}
\end{figure}

The interpolation problem of the two dimensional case is set up as tensor product of the one dimensional example as indicated in \myfigref{fig:exampleSchemes2d}.%
In particular, the basis and the definition of $\patchdomain$ are given by a tensor product of the univariate versions discussed in the previous section. 
The target function of the interpolation problem is given by
\begin{align}
	\label{eq:originalFunction2d}
	f(\uusurf,\vvsurf) &= \frac{1}{\sqrt{ \left(a_1-\uusurf\right)^2 + \left(a_2-\vvsurf\right)^2 }} && \with && a_1 = a_2 = -1.2 \:.
\end{align}
The knot insertion depth $\knotinsertion$ is set to $4$ for both intrinsic directions and the interpolation problem is performed for several trimming parameters $\trim$, $0.5 < \trim < 1$. 
For the untrimmed case $t=1$ the resulting $\kappa(\myMat{A}_{\uu})$ and $\myNorm{ \err_{rel}}_{\Ltwo}$ are summarized in~\mytabref{tab:untrimmedResults2d}.
All other results are shown in \myfigref{fig:example2d_cond} and \myfigref{fig:example2d_err}.%
\begin{mytableposition}
  {Results of the interpolation problem of the untrimmed bivariate basis.}%
  {untrimmedResults2d}
  {ccc}{htb}
  \mytableheader{
	$\pu$ & $\kappa$  & $\| \epsilon_{rel}\|_{L_2}$ }  
	$2$    & \num{6.25}         & \num{2.10822e-04} \\ 
	$3$    & \num{18.57448} & \num{4.48755e-05} \\ 
	$4$    & \num{63.01519} & \num{7.65653e-06} \\ 
\end{mytableposition}%

\tikzfigposition{tikz/resultsInterpolation2d_cond}{0.95}
{Condition number $\kappa(\myMat{A}_{\uu})$ of the bivariate basis for several degrees $\pu$ and trimming parameters~$\trim$. The labels of the horizontal axis indicate knots of the trimmed basis.}{fig:example2d_cond}{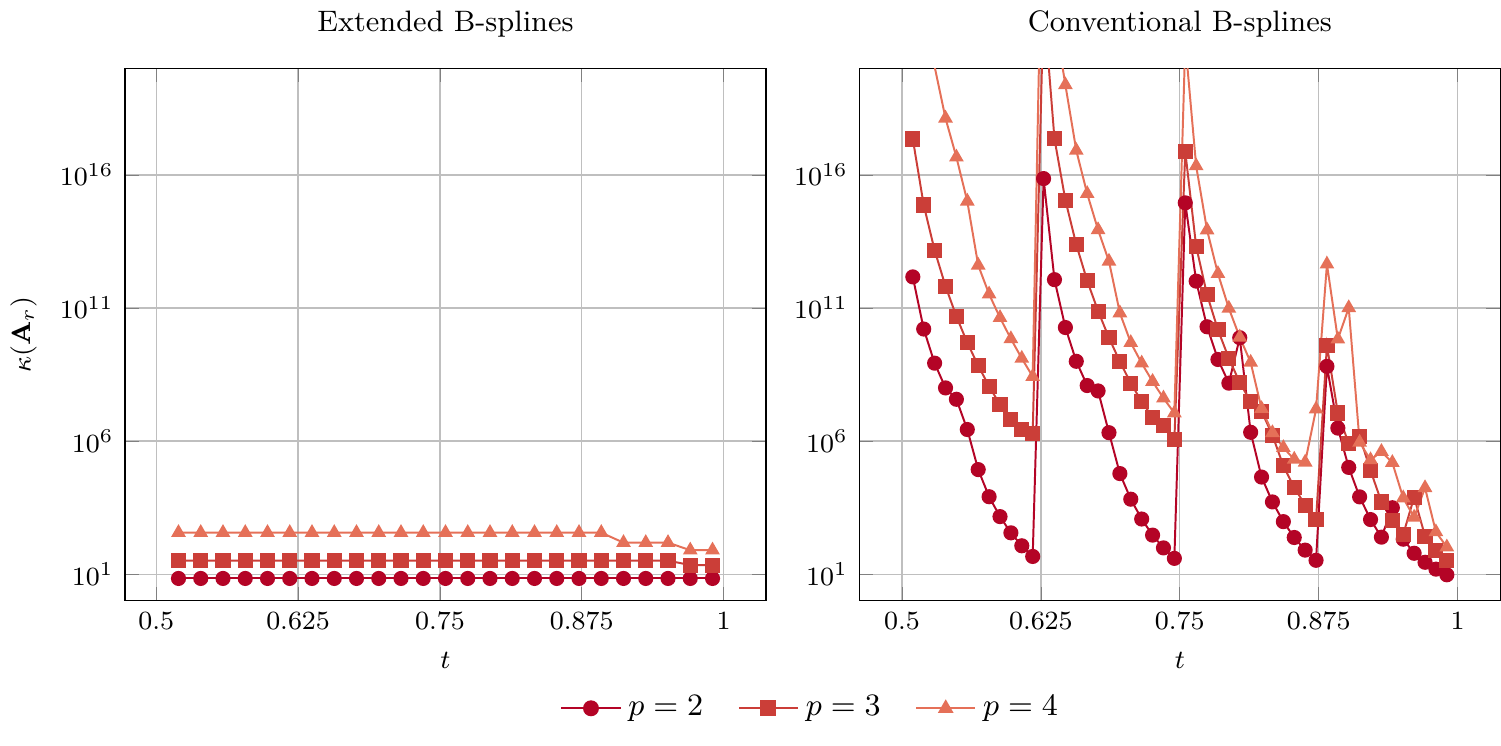}{th}%
\tikzfigposition{tikz/resultsInterpolation2d_err}{0.95}
{Relative interpolation error of the bivariate basis for several degrees~$\pu$ related to the trimming parameter~$\trim$. The labels of the horizontal axis indicate knots of the trimmed basis.}{fig:example2d_err}{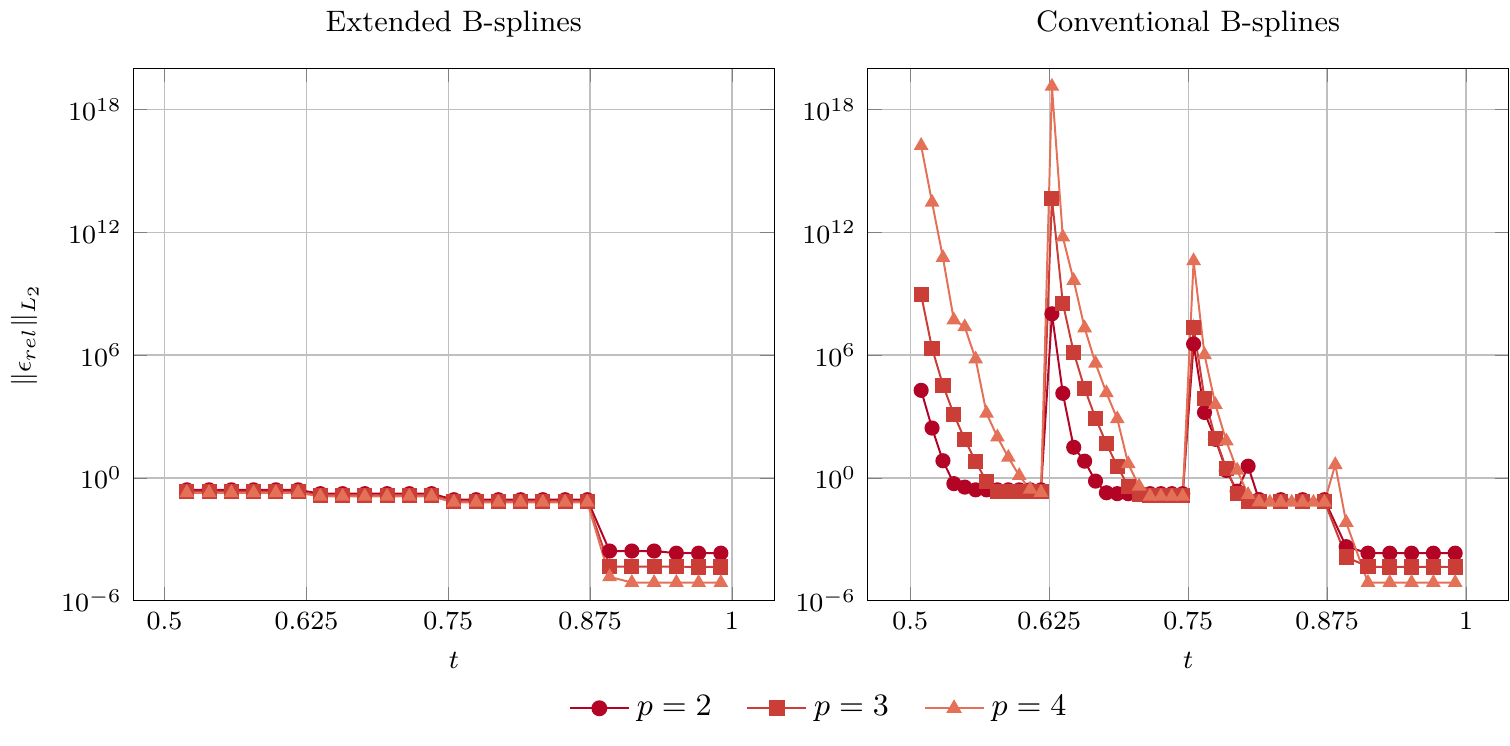}{h!}%

It is evident that the simple approach negatively affects the condition number and subsequently the quality of the approximation.
The error peaks near knot values are in fact disastrous.
Astonishingly, extended B-splines provide a stable basis where $\kappa(\myMat{A}_{\uu})$ barely changes for the same trimming cases. 
The reduction of the approximation accuracy occurs due to the reduction of the degrees of freedom~$\DOF$, i.e.~number of extended B-splines, as the trimming parameter $t \to 0.5$.

\subsection{Influence of Non-Uniformity}
\label{sec:nonUnifromTest}

\begin{figure}[h]
  \centering
  \begin{minipage}[b]{0.9\textwidth}
	\tikzsubfig{tikz/nonUniformExampleSetting2}{2.2}{$\trim=0.51$}{fig:example1dnonuniform_setting2}{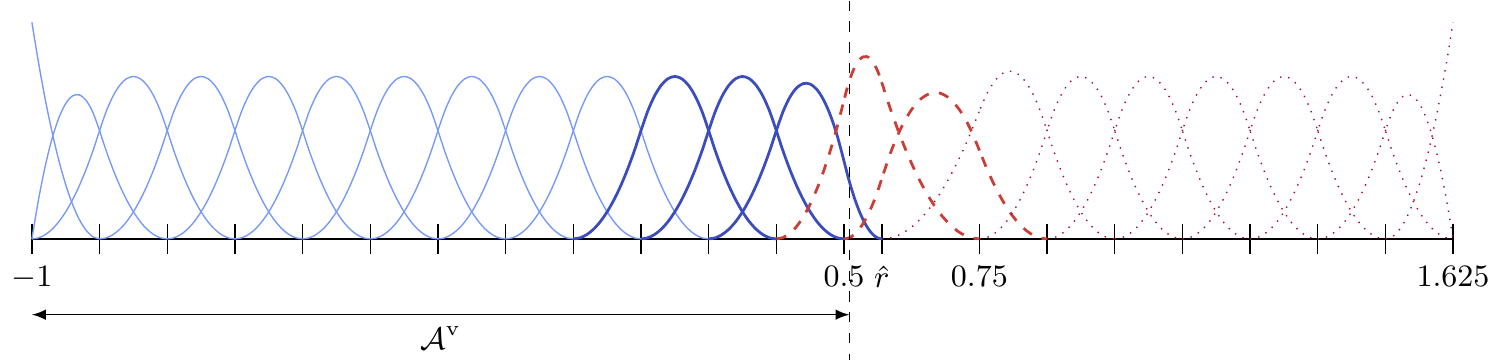}
  \end{minipage}
  \hfill
  \begin{minipage}[b]{0.9\textwidth}
	\tikzsubfig{tikz/nonUniformExampleSetting}{2.2}{$\trim=0.80$}{fig:example1dnonuniform_setting}{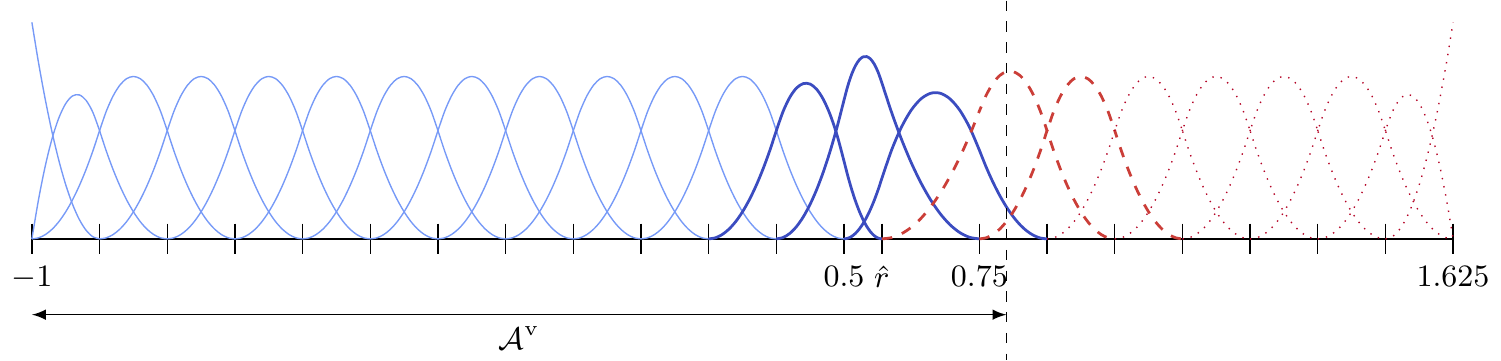}
  \end{minipage}
  \caption{
      Test setting to examine the influence of non-uniformity on the stabilization of extended B-splines: the size of (a)~the trimmed knot span or (b) the one that provides the stable B-splines is varied. In both cases, the size is determined by a variable knot $\variableKnot \in [ 0.5,0.75 ]$. Degenerated B-splines are indicated by dashed lines, while the thicker continuous lines represent stable B-splines used for the stabilization.
  }
  \label{fig:exampleSchemesNonUniform}
\end{figure}

In this study, the influence of the knot spacing to the extended B-spline stabilization is investigated.
Therefore, the following problem is considered:
an open knot vector defines a univariate parameter space from $-1.0$ to $1.625$ for various degrees $\pu=\{2,3,4\}$.
The interior knots are uniformly distributed with a knot span size of $h=0.125$.
This uniform pattern is interrupted by a single variable knot $\variableKnot \in [ 0.5,0.75 ]$.
Note that this set up includes the case of multiple knot values, in particular when $\variableKnot = 0.5$ or $\variableKnot = 0.75$.
Similar to the previous examples, a trimming parameter $\trim$ is introduced to specify the valid domain $\patchdomain \in [-1,\trim)$.
In particular, it is chosen in such a manner that either the trimmed knot span or the knot span $\indexSpan$ that provides the polynomial segments $\BsplineSeg^{\indexSpan}_{\indexA}$ for the stabilization is influenced by the variation of $\variableKnot$, i.e.~$\trim=0.51$ and $\trim=0.80$, respectively.
This is illustrated for a quadratic basis in \myfigref{fig:exampleSchemesNonUniform}.
The effectiveness of the stabilization is expressed by the condition number of the spline collocation matrix~$\kappa(\myMat{A}_\uu)$.
The results for both trimming cases are summarized in \myfigref{fig:example1dnonuniform_cond}.
The actual approximation error of this example is scarcely affected by the variation of $\variableKnot$ and is therefore omitted.
\tikzfig{tikz/resultsInterpolation1dnonuniform_cond}{0.95}
{Condition number $\kappa(\myMat{A}_{\uu})$ of the non-uniform basis for several degrees $\pu$ and two different trimming parameters~$\trim$. The horizontal axis indicates the position of the variable knot $\variableKnot$.}{fig:example1dnonuniform_cond}{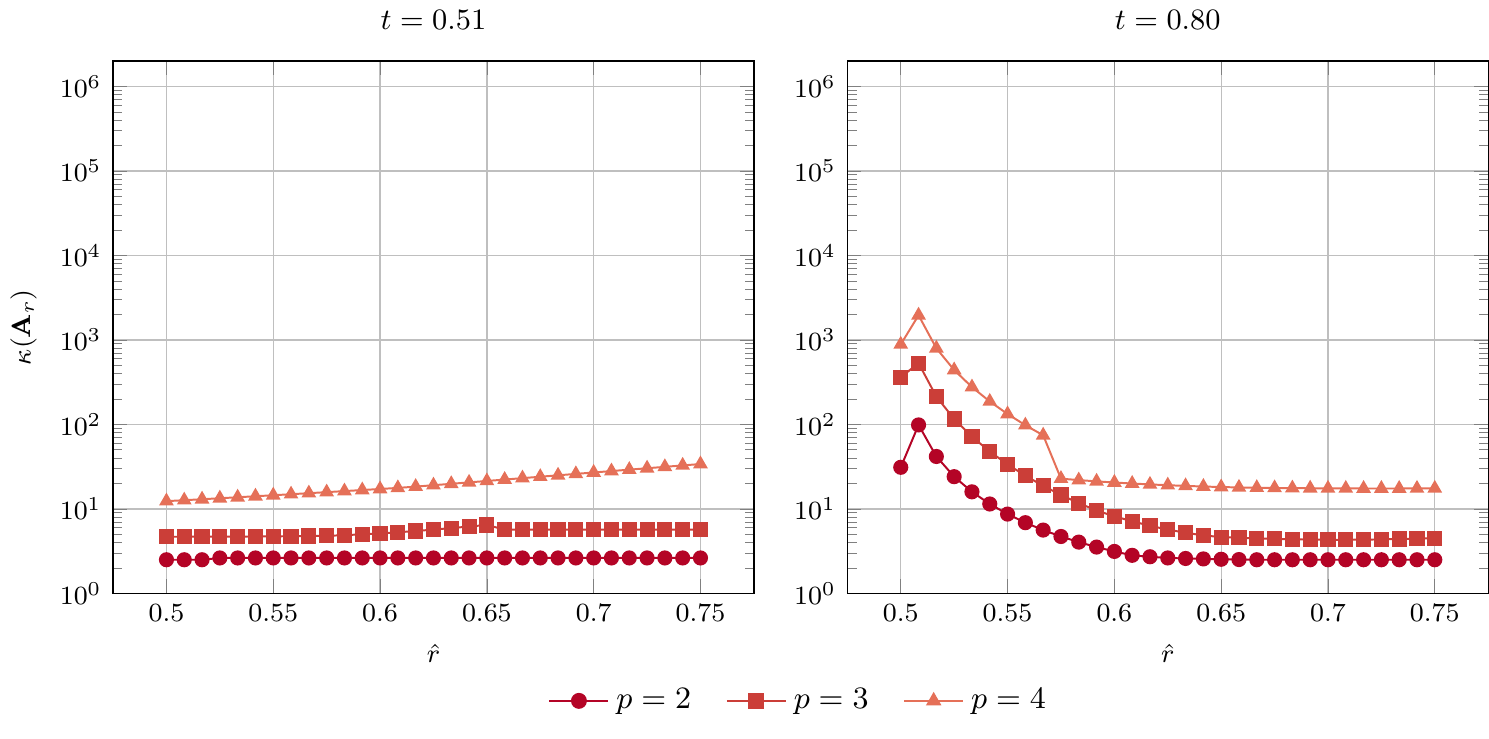}%

All condition numbers obtained are moderate, but 
it is apparent that $\kappa(\myMat{A}_{\uu})$ is affected by the knot spacing. 
Non-uniformity of the trimmed knot span barely changes $\kappa(\myMat{A}_{\uu})$ as shown by the trimming case $\trim=0.51$.
Looking at the case $\trim=0.80$ where the knot span $\indexSpan$ is determined by the variable knot, it can be concluded that knot spans should not become too small, if they provide the polynomial segments $\BsplineSeg^{\indexSpan}_{\indexA}$.
These results are in general agreement with the statement of \citet{Hoellig2003a} that the mesh-ratio influences the stabilization in case of finite spline spaces.

\tikzfig{tikz/resultsInterpolation1dnonuniform_cond_adaptive}{0.95}
{Condition number $\kappa(\myMat{A}_{\uu})$ of the non-uniform basis for several degrees $\pu$. The horizontal axis indicates the position of the variable knot $\variableKnot$. An adaptive refinement scheme is employed in order to obtain a better mesh-ratio.}{fig:example1dnonuniform_cond_adaptive}{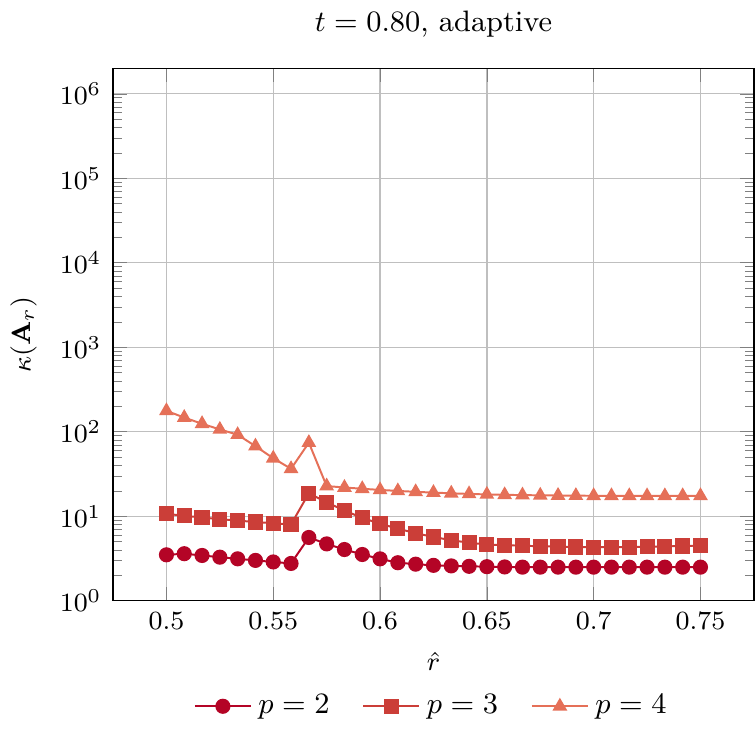}

In order to demonstrate that the mesh-ratio is important rather than the absolute size of the knot span~$\indexSpan$, another example is considered.
We repeat the study with $\trim=0.80$, but this time the adjacent knot span towards the trimming point is subdivided as soon as the ratio of these two knot spans is greater than $3$. 
The results of this adaptive scheme are shown in~\myfigref{fig:example1dnonuniform_cond_adaptive}.
The improvement due to the better mesh-ratio is clearly indicated by the kinks of the graphs at $\variableKnot\approx0.56$.
Finally, it should be emphasized that the stabilization is still independent of the trimming parameter $\trim$, i.e.~the variation of $\trim$ does not alter the results as long as it does not change the set of degenerated B-splines. 


\newpage
\subsection{Trimmed Cube}
\label{sec:stable}

\begin{figure}[tb]
  \centering
  \begin{minipage}[b]{0.45\textwidth}
	\tikzsubfig{tikz/3dcube_untrimmed.tex}{1.5}{Untrimmed}{fig:exampleCubeUntrimmed}{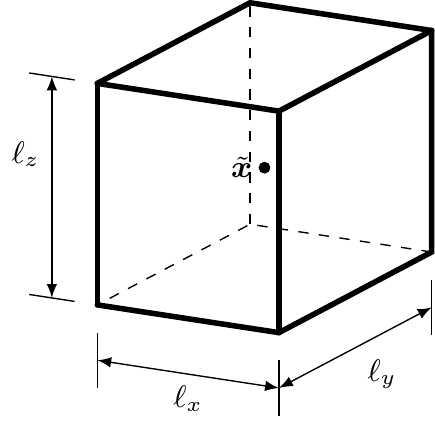}
  \end{minipage}
  \begin{minipage}[b]{0.45\textwidth}
	\tikzsubfig{tikz/3dcube_trimmed.tex}{1.5}{Trimmed}{fig:exampleCubeTrimmed}{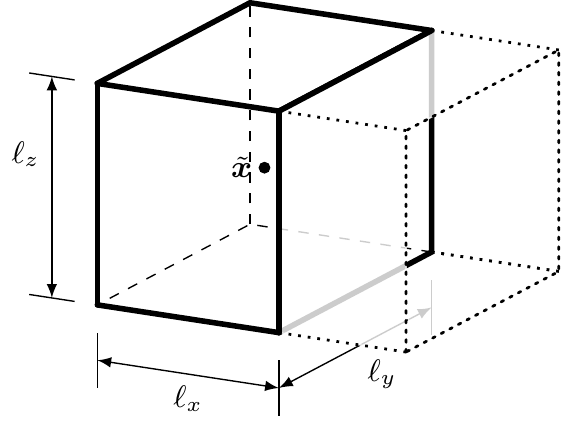}
  \end{minipage}
  \caption{Discretization of a unit cube by (a) regular patches and (b) trimmed patches.}
  \label{fig:exampleTrimmedCube}
\end{figure}

A unit cube is analyzed in order to investigate the approximation quality of extended B-splines in the context of an isogeometric BEM analysis.
The geometry is discretized by two different models as illustrated in~\myfigref{fig:exampleTrimmedCube}.%
One is described by \num{6} regular patches, whereas \num{4} trimmed patches are included in the other.
Both represent the same geometry, i.e.~$\ell_x=\ell_y=\ell_z=1.0$, which defines the boundary $\boundary$ of an infinite domain $\domain=\R^3 \setminus \domain^{-}$ with $\domain^{-}$ denoting the void.
The boundary condition is given by
\begin{align}
    \ttens(\boundarypt) & = \fund{T}(\sourcept,\boundarypt) &&
    \boundarypt\in\boundary \beistrich \sourcept\in\Omega^{-}.
\end{align}
In particular, a source point $\sourcept$ in the center of the cube defines the boundary conditions for the exterior Neumann problem. 
The Laplace as well as the \Lame-Navier equation is considered.
The discretizations are set up for different degrees $\pu=\left\{1,2,3\right\}$ and knot insertion is applied to improve the solutions.
The relative approximation error is determined by
\begin{align}
  \label{eq:direct-check}
  \err_{rel} & = \frac{\utens(\boundarypt) - \fund{U}(\tilde{\pt{x}},\pt{y})}{ \fund{U}(\tilde{\pt{x}},\pt{y}) } &&
    \forall\pt{y}\in\boundary \beistrich
    \tilde{\pt{x}}\in\Omega^{-}
\end{align}
where $\utens(\boundarypt)$ is the obtained solution.
The overall error is measured with respect to the \mbox{$\Ltwo$-norm}, i.e.~$\myNorm{\err_{rel}}_{\Ltwo}$.
The results are summarized for various degrees of freedom~$\DOF$ in \myfigref{fig:exampleTrimmedCube_err} and \myfigref{fig:exampleTrimmedCube_err_higherDegree}.

It can be observed from~\myfigref{fig:exampleTrimmedCube_err} that the trimmed model yields essentially the same results as in the untrimmed case, for degrees $\pu=\{1,2\}$.
The graphs related to degree $\pu=3$ show a similar convergence behavior, yet with a noticeable offset in favor of the untrimmed discretization.
This offset may occur due to the fact that the distance of the trimming curve to the inner knot span which provides the stable basis functions increases with the degree. 

\tikzfig{tikz/cube_exterior_trimmed_directNeumann_potential_discontinuous.tex}{0.95}{Relative $\Ltwo$-error of an exterior Neumann problem on the cube example with respect to the number of degrees of freedom $\DOF$.}{fig:exampleTrimmedCube_err}{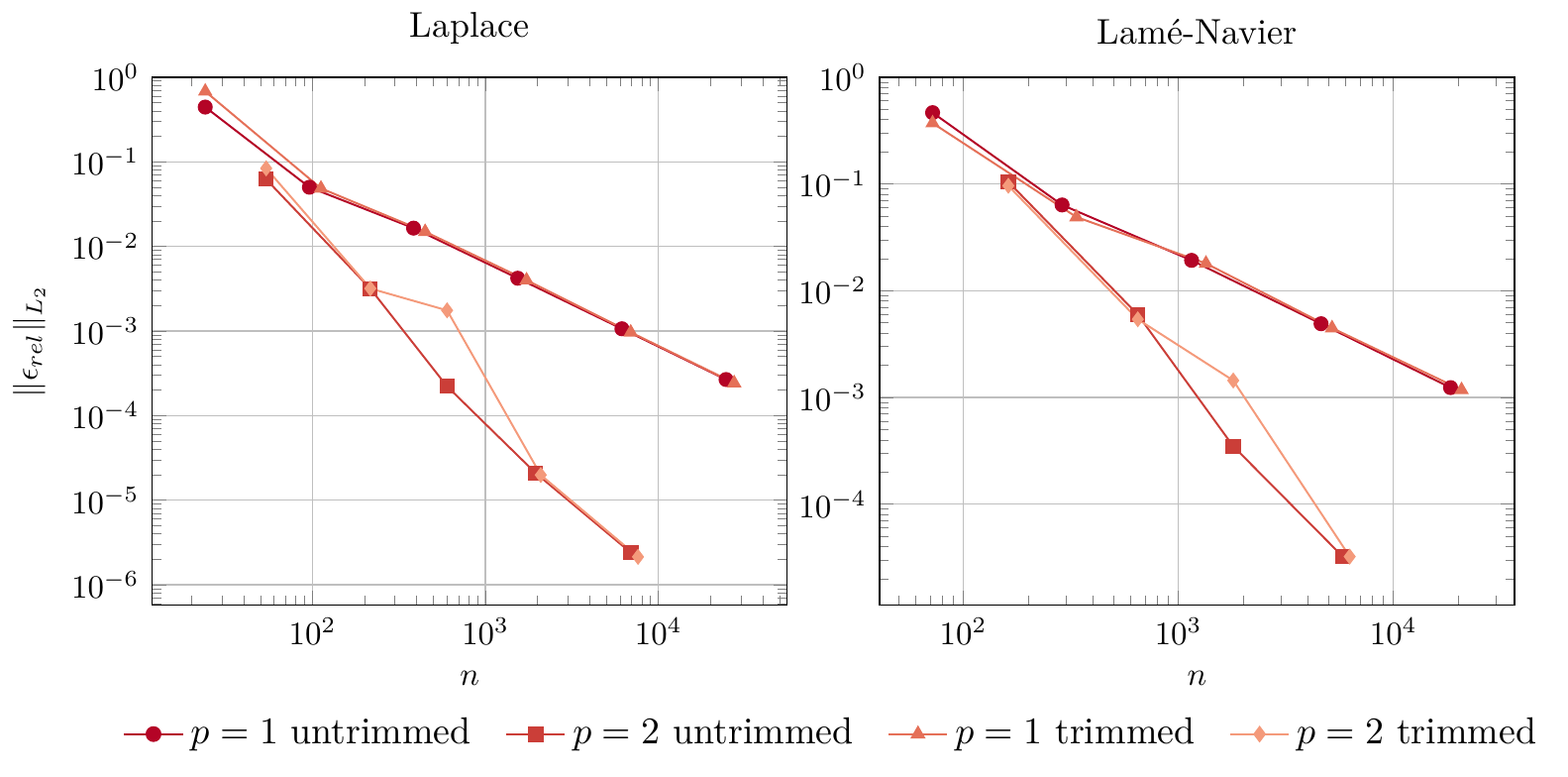}

\tikzfig{tikz/cube_exterior_trimmed_directNeumann_potential_discontinuous_higherDegree.tex}{0.95}{Relative $\Ltwo$-error of an exterior Neumann problem on the cube example discretized with basis functions of degree $\pu=3$.}{fig:exampleTrimmedCube_err_higherDegree}{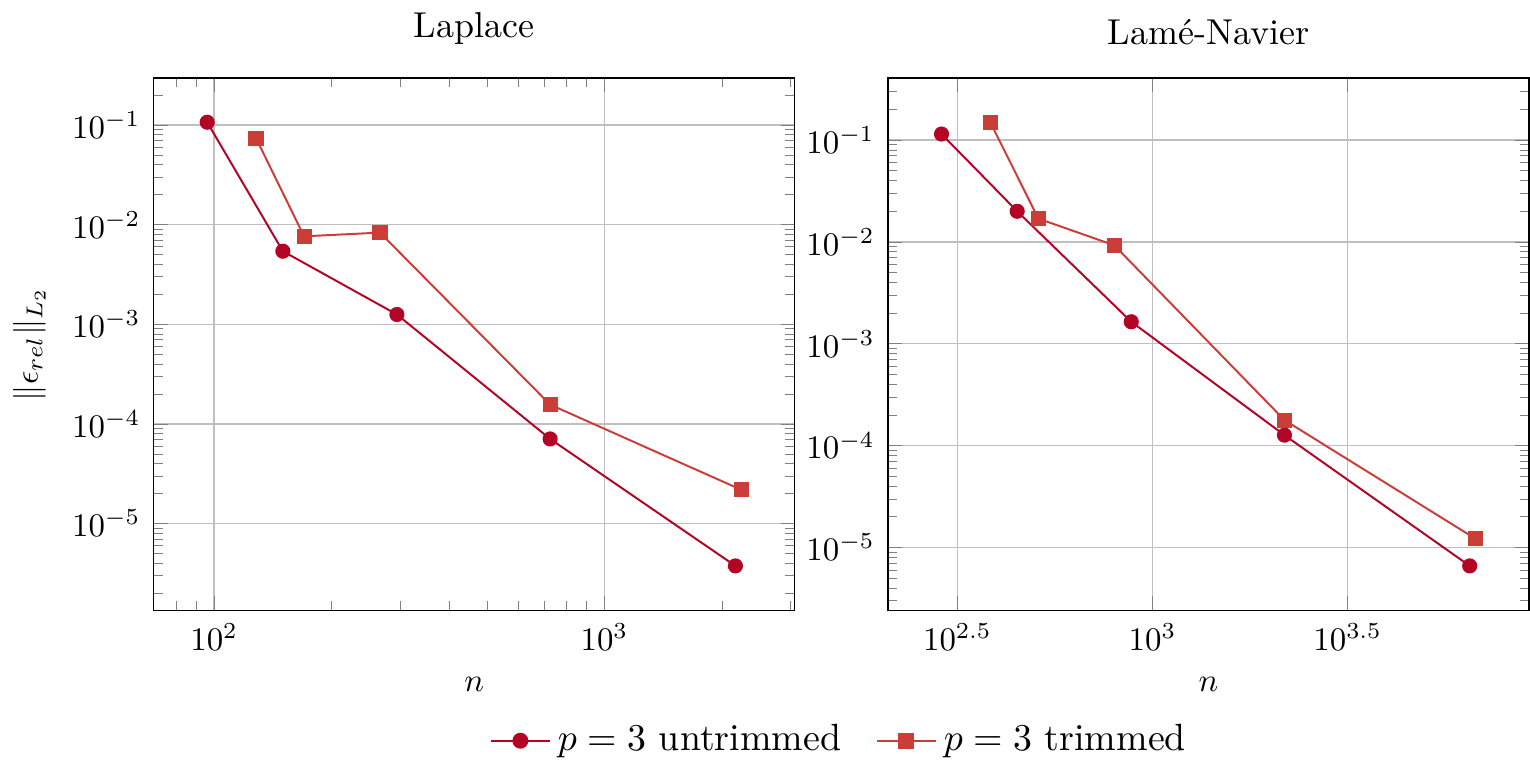}

\newpage
\subsection{Tunnel Cross Passage}
\label{sec:realworldExample}

\tikzfigposition{tikz/hudsonRiverTunnelModel.tex}{0.6}{CAGD model of the tunnel example with trimmed NURBS patches.}{fig:hudson_model}{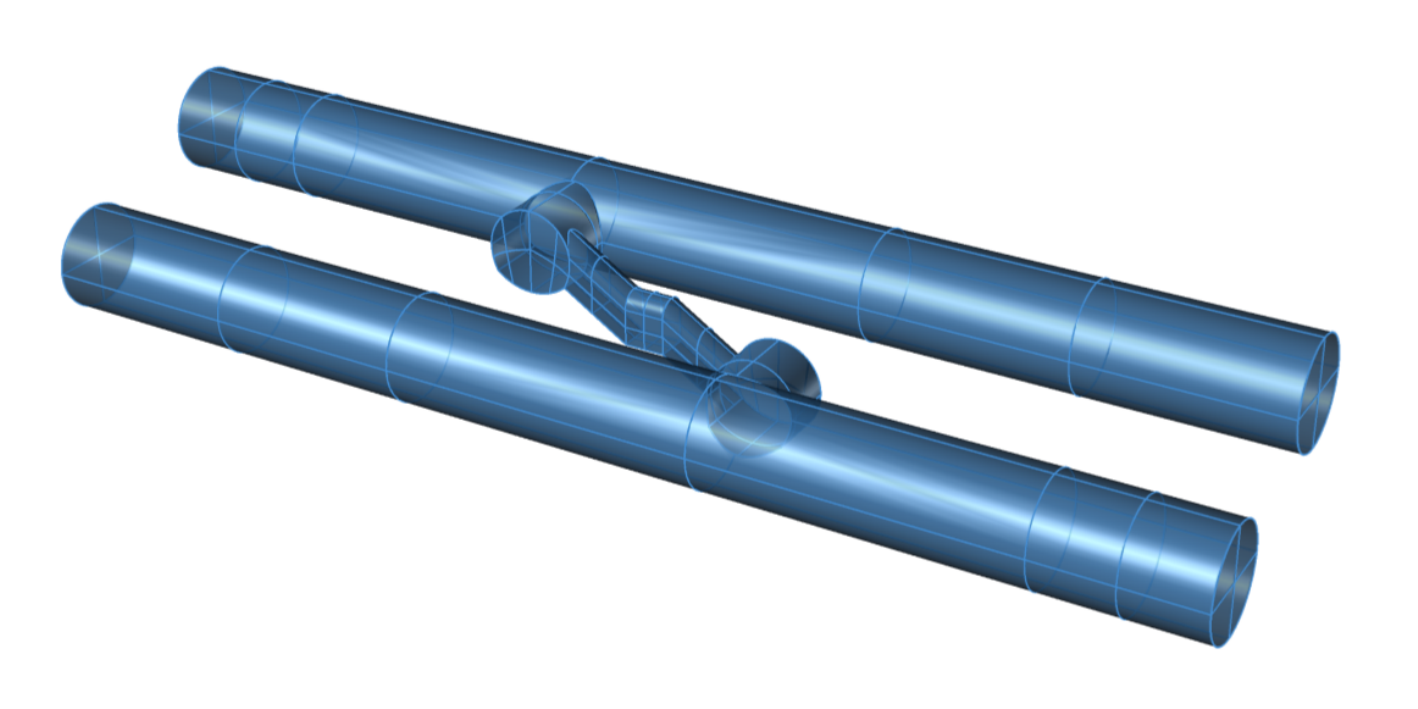}{b!}
In order to demonstrate the ability of the proposed scheme to handle complex trimming cases, the excavation of a metro tunnel is analyzed.
The geometry is specified by two parallel tunnel tubes which are connected by a cross passage as illustrated in \myfigref{fig:hudson_model}.
The CAGD model is defined by several trimmed  patches.
The tunnel tubes are cylinders, hence rational non-uniform B-splines (NURBS) provide the basis for the geometry description.
In order to apply extended B-splines for the analysis, 
two different approaches may be used: (i)~conversion of the CAGD model to a B-spline representation or (ii)~application of an independent field approximation~\cite{marussig2016b,Marussig2014aCMAME}.
In this paper, the latter is preferred since it allows to perform the analysis based on the original NURBS model without any geometrical approximations. 
The key idea of independent field approximation is to use different basis functions for the representation of the geometry and the approximations of the Cauchy data.
Hence, conventional B-splines are used for the discretization of the displacement and traction field.
This allows the straightforward application of extended B-splines.
In addition, the combination of NURBS for the geometry description and B-splines for the approximation of the Cauchy data has been shown to be more efficient \cite{marussig2016b} and does not lead to a loss of accuracy \cite{Li2011a,Marussig2014aCMAME}.
\begin{figure}[tb!]
  \centering
  \begin{minipage}[b]{0.32\textwidth}
\tikzsubfig{tikz/hudsonRiverTunnelModel_detailC.tex}{0.25}{Tunnel Tube}{fig:hudson_model_A}{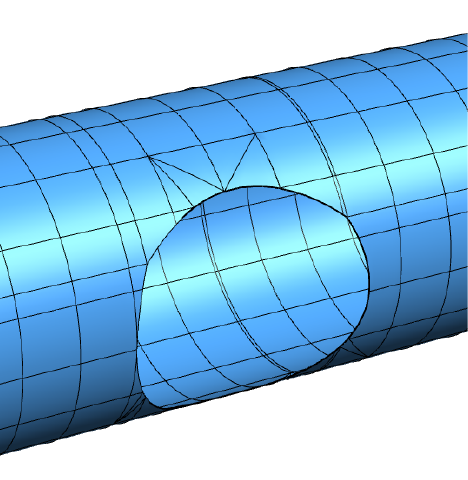}
  \end{minipage}
  \begin{minipage}[b]{0.32\textwidth}
\tikzsubfig{tikz/hudsonRiverTunnelModel_detailB.tex}{0.30}{Cross Passage}{fig:hudson_model_C}{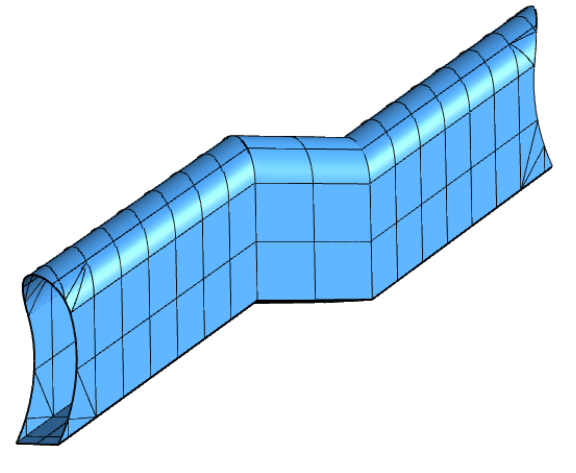}
  \end{minipage}
  \hfill
  \begin{minipage}[b]{0.32\textwidth}
\tikzsubfig{tikz/hudsonRiverTunnelModel_detailA.tex}{0.20}{Access Link}{fig:hudson_model_B}{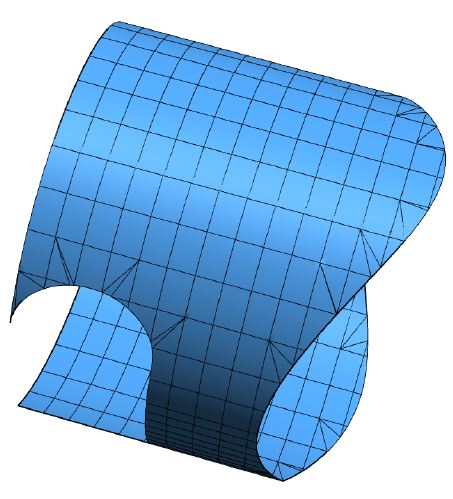}
  \end{minipage}
  \caption{Discretizations of different tunnel parts. Black lines indicate the applied partition into integration elements.}
  \label{fig:hudsonRiverDetails}
\end{figure}

In \myfigref{fig:hudsonRiverDetails}, closeups of some model parts together with the corresponding partition into integration elements are depicted.  
It should be noted that different trimming cases are covered including holes inside the patch, multiple trimming curves, and the representation of sharp features.

Linear elastic material properties are considered with
Young's modulus $E=\SI{313}{\mega\pascal}$ and Poisson ratio $\nu=0.2$.
The deformation of the tunnel due to its excavation in a single step is simulated. 
The boundary conditions are determined by applying excavation tractions related to the virgin stress field $\sigma_{xx}=\sigma_{yy}=\SI{1.375}{\mega\pascal}$ and $\sigma_{zz}=\SI{2.75}{\mega\pascal}$.
The known tractions along the boundary are calculated by multiplying the stress tensor with the outward normal of the surface.
The resulting displacements are visualized in~\myfigref{fig:hudson_result}.
This real-world example verifies the capability of extended B-splines to deal with complex trimming cases.

\tikzfig{tikz/hudsonRiverTunnelModel_result.tex}{0.18}{Resulting displacements of the tunnel example.}{fig:hudson_result}{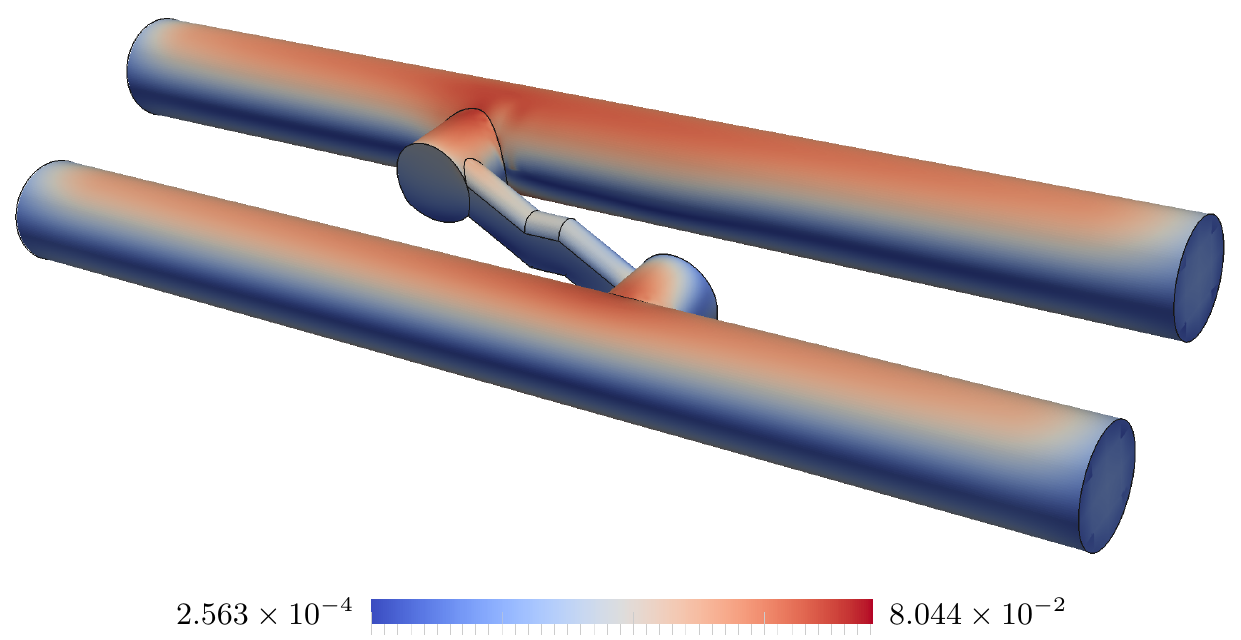}

  \section{Conclusion}
\label{sec:conclusion}

It is demonstrated that stability is a crucial aspect for the analysis of trimmed geometries.
Moreover, we examine the utilization of extended B-splines to resolve this issue.
These splines were originally introduced by Höllig~\cite{Hoellig2002a} and use linear combinations of conventional B-splines to establish a stable basis.
In this work, the concept of extended B-splines is derived from an engineering point of view rather than a mathematical one. 
The present procedure permits the straightforward construction of such splines for non-uniform parameter spaces and their application to collocation schemes.
The integration of extended B-splines into an existing isogeometric analysis software is simple, because it can be implemented as a supplementary task applied after the initially unstable system of equations has been set up.

The approximation quality and stabilization capability of the proposed approach have been verified by means of interpolation problems.
The stabilization is independent of the actual trimming position, established in a very flexible manner, and requires only the presence of a sufficient number of stable basis functions.
It is shown that the quality of the stabilization is negatively affected if knot spans that provide stable basis functions become too small in relation to their adjacent knot spans.
However, a proper knot span ratio can easily be achieved by refinement.
Additionally, extended B-splines have been studied
in the context of collocated isogeometric boundary element simulations.
The results obtained are very encouraging, in particular for lower degrees where trimmed discretizations yield the same quality of the results as untrimmed ones.  
For higher degrees, the convergence rate of the trimmed example is similar to the untrimmed reference, but with a considerable offset.

The concept of extended B-splines is very promising and we believe that it provides an essential feature in analyzing trimmed geometries.
The present scheme focuses on \mbox{B-splines}. 
Nevertheless, it is also applicable to geometries represented by NURBS since the stabilization affects the approximation of the physical fields only.
In particular, the analysis of NURBS models has been preformed using an independent field approximation, i.e.~the physical fields are described by extended B-splines while the geometrical mapping is still defined by conventional NURBS.
%
It is important to note that the stabilized B-spline basis proposed herein also applies to simulations using the finite element method.
In that case, the application of essential boundary conditions requires a more detailed discussion.

\section{Acknowledgment}

This research was supported by the Austrian Science Fund (FWF): P24974-N30. This support is gratefully acknowledged.

\appendix

  \section{Evaluation of Extrapolation Weights}
\label{chap:extrapolationWeightsExample}

\subsection{Explicit Representation}

The coefficients $\beta_\indexC$ of the polynomial $\Nbasis_{\indexB,\pu}$ 
are computed by \myeqref{eq:Nbasis_explicit_beta}. 
In particular, $\NSet_{\indexC,\indexG}$ is required which represents all $\indexC$-combinations with repetition of the knots $\left\{ \uu_{\indexB+1},\dots,\uu_{\indexB + \pu} \right\}$.
For example $\NSet_{\indexC,\indexG}$ in case of a cubic B-spline, i.e.~$\pu=3$, would be given by 
\begin{align*}
	\NSet_{3,1} &= \left\{ \uu_{\indexB + 1},\uu_{\indexB + 2},\uu_{\indexB + 3} \right\},
	&& \NSet_{2,1} = \left\{ \uu_{\indexB + 1},\uu_{\indexB + 2}  \right\}, 
	&& \NSet_{1,1} = \left\{ \uu_{\indexB + 1} \right\}, \\
	& &&\NSet_{2,2} = \left\{ \uu_{\indexB + 2},\uu_{\indexB + 3}  \right\}, 
	&&\NSet_{1,2} = \left\{ \uu_{\indexB + 2} \right\}, \\
	& &&\NSet_{2,3} = \left\{ \uu_{\indexB + 1},\uu_{\indexB + 3}  \right\}, 
	&&\NSet_{1,3} = \left\{ \uu_{\indexB + 3} \right\}.
\end{align*}

The power basis form \myRef{eq:Bspline_PowerBasis} of the polynomial segments~$\BsplineSeg^{\indexSpan}_{\indexA}$ is obtained by
\begin{align}
	\label{eq:Bspline_TaylorToPowerBasis}
	\BsplineSeg^{\indexSpan}_{\indexA}(\uu)
	= \sum^{\pu}_{\indexD=0} \alpha_\indexD\:(\uu-\taylorpoint)^{\indexD}
	= \poly_0 1+ \poly_1 \uu + \cdots + \poly_\pu \uu^\pu 
\end{align}	
where $\poly_\indexD$ is for $\pu = 1$:
\begin{align*}
	\poly_0 &=  1 \alpha_0 + 1 \alpha_1\left( -\taylorpoint \right), \\
	\poly_1 &=  1 \alpha_1, \\
	\shortintertext{\for \pu = 2:}  
	\poly_0 &=  1 \alpha_0 + 1 \alpha_1\left( -\taylorpoint \right) + 1 \alpha_{2} \left( -\taylorpoint \right)^{2},  \\
	\poly_1 &=  1 \alpha_1 + 2 \alpha_2 \left( -\taylorpoint \right), \\  
	\poly_2 &= 1 \alpha_2, \\ 
	\shortintertext{\for \pu = 3:}   
	\poly_0 &=  1 \alpha_0 + 1 \alpha_1\left( -\taylorpoint \right) + 1 \alpha_{2} \left( -\taylorpoint \right)^{2} + 1 \alpha_3 \left( -\taylorpoint \right)^3, \\
	\poly_1 &=  1 \alpha_1 + 2 \alpha_2 \left( -\taylorpoint \right) + 3 \alpha_{3} \left( -\taylorpoint \right)^{2},  \\  
	\poly_2 &= 1 \alpha_2 + 3 \alpha_3 \left( -\taylorpoint \right),  \\  
	\poly_3 &= 1 \alpha_3 \:.   	
\end{align*}
It should be noted that the factors in front of the coefficients $\alpha$ correspond to Pascal's triangle, which yields to the binomial coefficient in \myeqref{eq:Bspline_PowerBasis}.
In general, if the coordinate $\uu$ does not contribute to the derivative of a power basis form, e.g.~$\uu=0$, its evaluation simplifies to  
\begin{align}
	\label{eq:powerBasisDerivZero}
	\BsplineSeg^{\indexSpan^{\left( \indexC\right)}}_{\indexA}(0) & = \indexC! \: \poly_\indexC
\end{align}
which has been utilized to obtain \myeqref{eq:deBoorFix_explicit}.

\subsection{Example}

The example shown in~\myfigref{fig:extendedBsplineConcept} is considered to clarify the computation of extrapolation weights $\eweight_{{\indexA},{\indexB}}$.
The knot vector is given by $\KV = \left\{1,1,1,2,3,4,4,4\right\}$, the first B-spline is degenerated, i.e.~$\degSet = \{0\}$, and the knot span $\indexSpan = 3$ is the closest non-trimmed interval.
Hence, the polynomial $\Nbasis_{0,2}$ is determined by the knot values $\left\{1,1\right\}$ and the corresponding coefficients are $\beta_\indexC =\left\{1,-2,1\right\}$.
The polynomial segments $\BsplineSeg^{\indexSpan}_{\indexA}$ obtained by \myeqref{eq:Bspline_PowerBasis} are
\begin{align}
\label{eq:segment1}
\BsplineSeg^{3}_{1} &= 0.5\:\uu^2 - 3\: \uu + 4.5, \\
\label{eq:segment2}
\BsplineSeg^{3}_{2} &= -1\: \uu^2 + 5\: \uu - 5.5, \\
\label{eq:segment3}
\BsplineSeg^{3}_{3} &= 0.5\: \uu^2 - 2\: \uu + 2 \: .
\end{align}
These polynomials are the target functions for an interpolation problem and the corresponding results provide the extrapolation weights~$\eweight_{{\indexA},{\indexB}}$.
In the following, this is demonstrated by three different approaches:
\begin{itemize}
    \item Spline interpolation
    \item Direct evaluation of the functional~\myRef{eq:EW} by applying Horner's method to the explicit representations~\myRef{eq:Nbasis_explicit} and~\myRef{eq:Bspline_PowerBasis} of the polynomials $\Nbasis_{\indexB,\pu}$ and $\BsplineSeg^{\indexSpan}_{\indexA}$
    \item Indirect evaluation of the functional using the simplified expression~\myRef{eq:deBoorFix_explicit}
\end{itemize}
The spline interpolation approach is included since it is perhaps more familiar than quasi interpolation. In particular, it shall emphasize two points: (i) extrapolation weights are indeed obtained by solving an interpolation problem and (ii) extrapolation weights of non-degenerated B-splines are trivial, i.e.~either $0$ or $1$, as discussed in \mysecref{sec:ExtendedBsplines1Dweights}.    
In general, quasi interpolation should be preferred to obtain $\eweight_{{\indexA},{\indexB}}$ since it allows to compute solely the non-trivial results.
Moreover, the anchors~$\abscissa$ needed for the spline interpolation have to be within the trimmed knot span, which is not optimal.
The \emph{indirect} evaluation of the functional is a simplification of the \emph{direct} one, because it requires  only the coefficients of the polynomials $\Nbasis_{\indexB,\pu}$ and $\BsplineSeg^{\indexSpan}_{\indexA}$ rather than their evaluation at a point~$\lambdapoint_\indexB$.

\paragraph{Spline Interpolation}

The anchors are chosen to be $\abscissa=\left\{1,1.5,2\right\}$ leading to the spline collocation matrix $\myMat{A}_\uu$ and the right hand side $\myVec{f}_\indexA$ for each \myRef{eq:segment1} -- \myRef{eq:segment3}:
\begin{align*}
\myMat{A}_\uu &=
	\begin{pmatrix}
		1    & 0      & 0      \\            
		0.25 & 0.625  & 0.125  \\            
		0    & 0.5    & 0.5                
	\end{pmatrix},
	&& \myVec{f}_1 = 
	\begin{pmatrix}
		2 \\ 1.125 \\ 0.5          
	\end{pmatrix},
	&& \myVec{f}_2 = 
	\begin{pmatrix}
		-1.5 \\ -0.25 \\ 0.5          
	\end{pmatrix},
	&& \myVec{f}_3 = 
	\begin{pmatrix}
		0.5 \\ 0.125 \\ 0            
	\end{pmatrix}.
\end{align*}
Solving the system of equations gives 
\begin{align*}
\myMat{M} &=
	\begin{pmatrix}
		\eweight_{1,0}  & \eweight_{1,1} & \eweight_{1,2}    \\            
		\eweight_{2,0}  & \eweight_{2,1} & \eweight_{2,2}    \\            
		\eweight_{3,0}  & \eweight_{3,1} & \eweight_{3,2}              
	\end{pmatrix}
	=
	\begin{pmatrix}
		2  & 1 & 0      \\            
	     -1.5  & 0 & 1  \\            
	      0.5  & 0 & 0                
	\end{pmatrix}
\end{align*}
where the matrix-rows correspond to \myRef{eq:segment1} -- \myRef{eq:segment3} and the first column provides the sought extrapolation weights $\eweight_{{\indexA},{\indexB}}$ for $\degSet = \{0\}$.

\paragraph{Direct Evaluation of the Functional}

The explicit representation of $\Nbasis_{\indexB,\pu}$ can be evaluated by Horner's method for higher derivatives \cite{www:horner}.
The procedure is summarized in~\myalgref{alg:horner}, which returns all values $ \Nbasis^{\left(\indexD\right)}_{\indexB,\pu}$ with $\indexD = \left\{0,\dots,\pu \right\}$ collected in a vector $\myVec{v}^\Nbasis$.
The same procedure can be applied to the explicit representation of $\BsplineSeg^{\indexSpan}_{\indexA}$, if $\beta_\indexC$ is substituted by $\poly_\indexC$.
For the given example the position of the evaluation point is chosen to be $\lambdapoint_\indexB=1$ which leads to the following values
\begin{align*}
 \myVec{v}^{\BsplineSeg}_1 &=  \begin{pmatrix} 2 \\ -2 \\ 1 \end{pmatrix},
 && \myVec{v}^{\BsplineSeg}_2 =  \begin{pmatrix} -1.5 \\ 3 \\ -2 \end{pmatrix},
 && \myVec{v}^{\BsplineSeg}_3 =  \begin{pmatrix} 0.5 \\ -1 \\ 1 \end{pmatrix}
 && \und && 
 \myVec{v}^\Nbasis =  	\begin{pmatrix} 0 \\ 0 \\ 2 \end{pmatrix}.
\end{align*}
Applying these values to \myRef{eq:EW} provides the extrapolation weights
\begin{align*}
	\eweight_{1,0}  &= \frac{1}{2} \left[ 2 \cdot 2  \right] = 2,
	&& \eweight_{2,0}  = \frac{1}{2} \left[ 2 \cdot (-1.5)  \right] = -1.5,
	&& \eweight_{3,0}  = \frac{1}{2} \left[ 2 \cdot 0.5  \right] = 0.5 \:.
\end{align*}

\begin{myalgorithm}{Horner's method for higher derivatives}{alg:horner}
	\REQUIRE polynomial coefficients $\beta_\indexC$ of $\Nbasis_{\indexB,\pu}$ and the coordinate $\lambdapoint_\indexB$
	\STATE initialize matrix $\myMat{M} \in \R^{\pu+2 \times \pu+1}$ 
	\FOR{$\indexC=0$ \TO $\pu$}
	\STATE $\myMat{M}_{0,\indexC} \store \beta_\indexC$
	\ENDFOR
	\FOR[compute Horner coefficients]{$\indexD=0$ \TO $\pu$} 
	\STATE $\myMat{M}_{\indexD+1,\pu} = \myMat{M}_{\indexD,\pu}$
	\STATE $\indexE \store \pu-1$
	\WHILE{$\indexE \geq \indexD$}
	\STATE $\myMat{M}_{\indexD+1,\indexE} = \lambdapoint_\indexB \: \myMat{M}_{\indexD+1,\indexE+1} + \myMat{M}_{\indexD,\indexE}$ 
	\STATE $\indexE \store \indexE - 1$
	\ENDWHILE
	\ENDFOR
	\STATE initialize vector $\myVec{v}^{\Nbasis} \in \R^{\pu+1}$
	\FOR[compute derivatives]{$\indexD=0$ \TO $\pu$ }
	\STATE $ \myVec{v}^\Nbasis_\indexD  \store \indexD! \: \myMat{M}_{\indexD+1,\indexD} $
	\ENDFOR
	\RETURN $\myVec{v}^\Nbasis $
\end{myalgorithm}%

\paragraph{Indirect Evaluation of the Functional} 
For the given example the coefficients  $\Nbasis^{\left(\pu-\indexC\right)}_{\indexB,\pu} = \left( \pu - \indexC \right)! \: \beta_{\pu-\indexC}$ and $\BsplineSeg^{\indexSpan^{\left( \indexC\right)}}_{\indexA} = \indexC ! \:  \poly_\indexC$ of \myeqref{eq:deBoorFix_explicit} are
\begin{align*}
	\Nbasis_{0,2} &= 1 \cdot 1, &&
	\Nbasis^{\left(1\right)}_{0,2} = 1 \cdot (-2), && 
	\Nbasis^{\left(2\right)}_{0,2} = 2 \cdot 1,
	\shortintertext{respectively} 
	\BsplineSeg^{3}_{1} &= 1 \cdot 4.5 ,&&
	\BsplineSeg^{3^{\left( 1\right)}}_{1} = 1 \cdot (-3) ,&&
	\BsplineSeg^{3^{\left( 2\right)}}_{1} = 2 \cdot 0.5,      \\
	\BsplineSeg^{3}_{2} &= 1 \cdot (-5.5) ,&&
	\BsplineSeg^{3^{\left( 1\right)}}_{2} = 1 \cdot 5 ,&&
	\BsplineSeg^{3^{\left( 2\right)}}_{2} = 2 \cdot  (-1),    \\
	\BsplineSeg^{3}_{3} &= 1 \cdot 2 ,&&
	\BsplineSeg^{3^{\left( 1\right)}}_{3} = 1 \cdot (-2),&&
	\BsplineSeg^{3^{\left( 2\right)}}_{3} = 2 \cdot  0.5  \: .  
\end{align*}
Hence, the extrapolation weights are computed by
\begin{align*}
	\eweight_{1,0}  &= \frac{1}{2} \left[  2 \cdot 4.5 - (-2) \cdot (-3) + 1 \cdot 1 \right] = 2, \\
	\eweight_{2,0}  &= \frac{1}{2} \left[ 2 \cdot (-5.5) - (-2) \cdot 5 + 1 \cdot (-2)  \right] = -1.5, \\ 
	\eweight_{3,0}  &= \frac{1}{2} \left[ 2 \cdot 2 - (-2) \cdot (-2)  + 1 \cdot 1  \right] = 0.5 \: .
\end{align*}

  \bibliographystyle{plainnat}

\bibliography{localbib}




\end{document}